\pdfoutput=1
\RequirePackage{ifpdf}
\ifpdf 
\documentclass[pdftex]{sigma}
\else
\documentclass{sigma}
\fi

\usepackage{fancyhdr}
\usepackage[matrix,arrow,tips]{xy}
\usepackage{enumerate}
\newtheorem{Theorem}{Theorem}[section]
\newtheorem{Corollary}[Theorem]{Corollary}
\newtheorem{Conjecture}[Theorem]{Conjecture}
\newtheorem{Lemma}[Theorem]{Lemma}
\newtheorem{Proposition}[Theorem]{Proposition}
{ \theoremstyle{definition}
\newtheorem{Definition}[Theorem]{Definition}
\newtheorem{Remark}[Theorem]{Remark} }
\numberwithin{equation}{section}
\DeclareMathOperator{\End}{End}
\DeclareMathOperator{\tr}{tr}

\begin{document}

\newcommand{\tHe}[3]{\,_{2}H_{1} \left(\genfrac{.}{.}{0pt}{}{#1}{#2};#3\right)}

\allowdisplaybreaks

\newcommand{\arXivNumber}{1312.6577}

\renewcommand{\PaperNumber}{113}

\FirstPageHeading

\ShortArticleName{Matrix Valued Classical Pairs Related to Compact Gelfand Pairs of Rank One}

\ArticleName{Matrix Valued Classical Pairs\\
Related to Compact Gelfand Pairs of Rank One}

\Author{Maarten VAN PRUIJSSEN~$^\dag$ and Pablo ROM\'AN~$^\ddag$}

\AuthorNameForHeading{M.~van Pruijssen and P.~Rom\'an}

\Address{$^\dag$~Universit\"at Paderborn, Institut f\"ur Mathematik,\\
\hphantom{$^\dag$}~Warburger Str.~100, 33098 Paderborn, Germany}
\EmailD{\href{mailto:vanpruijssen@math.upb.de}{vanpruijssen@math.upb.de}}
\URLaddressD{\url{http://www.mvanpruijssen.nl}}

\Address{$^\ddag$~CIEM, FaMAF, Universidad Nacional de C\'ordoba, Medina Allende s/n Ciudad Universitaria,\\
\hphantom{$^\ddag$}~C\'ordoba, Argentina}
\EmailD{\href{mailto:roman@famaf.unc.edu.ar}{roman@famaf.unc.edu.ar}}
\URLaddressD{\url{http://www.famaf.unc.edu.ar/~roman}}

\ArticleDates{Received April 30, 2014, in f\/inal form December 12, 2014; Published online December 20, 2014}

\Abstract{We present a~method to obtain inf\/initely many examples of pairs $(W,D)$ consis\-ting of a~matrix weight~$W$ in
one variable and a~symmetric second-order dif\/ferential ope\-ra\-tor~$D$.
The method is based on a~uniform construction of matrix valued polynomials starting from compact Gelfand pairs~$(G,K)$
of rank one and a~suitable irreducible~$K$-representation.
The heart of the construction is the existence of a~suitable base change~$\Psi_{0}$.
We analyze the base change and derive several properties.
The most important one is that $\Psi_{0}$ satisf\/ies a~f\/irst-order dif\/ferential equation which enables us to compute the
radial part of the Casimir operator of the group~$G$ as soon as we have an explicit expression for~$\Psi_{0}$.
The weight~$W$ is also determined by~$\Psi_{0}$.
We provide an algorithm to calculate~$\Psi_{0}$ explicitly.
For the pair $(\mathrm{USp}(2n),\mathrm{USp}(2n-2)\times\mathrm{USp}(2))$ we have implemented the algorithm in GAP so
that individual pairs $(W,D)$ can be calculated explicitly.
Finally we classify the Gelfand pairs $(G,K)$ and the~$K$-representations that yield pairs $(W,D)$ of size $2\times2$
and we provide explicit expressions for most of these cases.}


\Classification{22E46; 33C47}

\vspace{-2mm}

\section{Introduction}

Matrix valued orthogonal polynomials (MVOPs) in one variable are generalizations of scalar valued orthogonal polynomials
and they already show up in the 1940s~\cite{Krein1,Krein2}.
Since then, MVOPs have been studied in their own right and they have been applied and studied in dif\/ferent f\/ields such
as scattering theory, spectral analysis and representation theory~\cite{Berez,Damanik,
Geronimo,Groenevelt-Ismail-Koelink,Groenevelt-Koelink}.
In this paper we are concerned with obtaining families of MVOPs whose members are simultaneous eigenfunctions of
a~symmetric second-order dif\/ferential operator.

\looseness=-1
Fix $N\ge1$ and an interval $I\subset\mathbb{R}$.
We write $\mathbb{M}=\End(\mathbb{C}^{N})$ and the Hermitian adjoint of $A\in\mathbb{M}$ is denoted by $A^{*}$.
The space $\mathbb{M}[x]$ is an $\mathbb{M}$-bimodule.
A~matrix weight is a~function $W:I\to\mathbb{M}$ with f\/inite moments and $W(x)=W(x)^*$ and $W(x)>0$ almost everywhere.
The pairing
\begin{gather*}
\mathbb{M}[x]\times\mathbb{M}[x]\to\mathbb{M}: \ (P,Q)\mapsto\langle P,Q\rangle_{W}:=\int_{I}P(x)^{*}W(x)Q(x)dx
\end{gather*}
is an $\mathbb{M}$-valued inner product that makes $\mathbb{M}[x]$ into a~right pre-Hilbert $\mathbb{M}$-module.
A~family of MVOPs (with respect to this pairing) is a~family $(P_{n}:\, n\in\mathbb{N})$ with $P_{n}\in\mathbb{M}[x]$
satis\-fying (1)~$\deg(P_{n})=n$, (2)~the leading coef\/f\/icient of $P_{n}$ is invertible, for all $n\in\mathbb{N}$ and (3)~$\langle P_{n},P_{m}\rangle_{W}=M_{n}\delta_{m,n}$ for all $n,m\in\mathbb{N}$.
Existence of such a~family is guaranteed by application of the Gram--Schmidt process on $(1,x,x^{2},\ldots)$.
Moreover, up to right multiplication by $\mathrm{GL}(\mathbb{C}^{N})$, a~family of MVOPs is uniquely determined by the
weight~$W$.

In~\cite{Duran97} the question was raised if there exists a~matrix weight together with a~dif\/ferential operator of
degree two that has the corresponding family of orthogonal polynomials as a~family of simultaneous eigenfunctions.
If $N=1$ then the answer is well known, we get the classical orthogonal polynomials~\cite{Bochner}.
In fact, the algebra of dif\/ferential operators that have the classical polynomials as simultaneous eigenfunctions is
a~polynomial algebra generated by a~second-order dif\/ferential operator, see e.g.~\cite{Miranian}.

In general, the algebra of dif\/ferential operators that act on the $\mathbb{M}$-valued polynomials is
$\End(\mathbb{M})[x,\partial_{x}]$.
We identify $\End(\mathbb{M})=\mathbb{M}\otimes\mathbb{M}$ such that a~simple tensor $A\otimes B$ acts on an
element $C\in\mathbb{M}$ via $(A\otimes B)C=ACB^{*}$.
A~polynomial $P\in\mathbb{M}[x]$ is an eigenfunction of a~dif\/ferential operator
$D\in\End(\mathbb{M})[x,\partial_{x}]$ if there exists an element $\Lambda\in\mathbb{M}$ such that
$DP=P\Lambda$.
This is justif\/ied by the fact that we consider $\mathbb{M}[x]$ as right pre-Hilbert $\mathbb{M}$-module.

A pair $(W,D)$ consisting of a~matrix weight and an element $D\in\End(\mathbb{M})[x,\partial_{x}]$ of order two
that is symmetric with respect to $\langle\cdot,\cdot\rangle_{W}$ is called a~matrix valued classical pair (MVCP).
Any family of MVOPs is automatically a~family of simultaneous eigenfunctions.

Consider the two subalgebras $(\mathbb{M}\otimes\mathbb{C})[x,\partial_{x}]$ and
$(\mathbb{C}\otimes\mathbb{M})[x,\partial_{x}]$ of $\End(\mathbb{M})[x,\partial_{x}]$.
If $(W,D)$ is a~classical pair with $D\in(\mathbb{C}\otimes\mathbb{M})[x,\partial_{x}]$ then the weight can be
diagonalized by a~constant matrix~\cite{Duran97}.
After publication of~\cite{Duran97} examples of MVCPs $(W,D)$ with $D\in(\mathbb{M}\otimes \mathbb{C})[x,\partial_{x}]$
came about, arising from analysis on compact homogeneous spaces in a~series of papers starting in~\cite{GPT} and ending
in~\cite{PT}.

A uniform construction of MVCPs arising from the representation theory of compact Lie groups is presented
in~\cite{HvP,vP} and was inspired by~\cite{KvPR,KvPR2,Koornwinder,Vretare}.
The input datum is a~compact Gelfand pair $(G,K)$ of rank one and a~certain face~$F$ of the Weyl chamber of~$K$.
For each $\mu\in F$, the output is an orthogonal family of $\mathbb{M}^{\mu}=\End(\mathbb{C}^{N_{\mu}})$ valued
functions $(\Psi_{d}^{\mu}:\, d\in\mathbb{N})$ on the circle $S^{1}$ and a~commutative algebra of dif\/ferential operators
$\mathbb{D}(\mu)$ for which the functions $\Psi_{d}^{\mu}$ are simultaneous eigenfunctions.
Here, $N_{\mu}$ is a~natural number that depends on the weight $\mu\in F$.
Moreover the functions $\Psi_{d}^\mu$ are determined by this property and a~normalization.
We call $\Psi_{d}^\mu$ the {\it full spherical function} of type~$\mu$ and degree~$d$.
The datum $(G,K,F)$ for which this construction applies is called a~multiplicity free system and they are classif\/ied in~\cite{HvP}.

We obtain families of MVOPs by multiplying the functions $\Psi^{\mu}_{d}$ from the left with the inverse of
$\Psi^{\mu}_{0}$.
The matrix weight $W^{\mu}$ is expressed in terms of $\Psi_{0}^{\mu}$, some data from the irreducible~$K$-representation
and a~scalar Jacobi weight that is associated to the symmetric space $G/K$.
Conjugating the elements of $\mathbb{D}(\mu)$ with $\Psi^{\mu}_{0}$ yields a~commutative algebra $\mathbb{D}^{\mu}$ of
dif\/ferential operators for which the MVOPs are simultaneous eigenfunctions.
In fact, the MVOPs are determined by this property and a~normalization.
The commutative algebra $\mathbb{D}^{\mu}$ is contained in $(\mathbb{M}\otimes \mathbb{C})[x,\partial_{x}]$.

The exact relation between $W^{\mu}$ and $\mathbb{D}^{\mu}$ is not yet understood on the level of the polynomials,
i.e.~it is not clear what exactly characterizes $\mathbb{D}^{\mu}$.
From this point of view it is not clear what the right notion of a~matrix valued classical pair should be.
In the spirit of classical orthogonal polynomials the MVOPs should be determined as eigenfunctions of a~commutative
algebra of dif\/ferential operators.
The algebras $\mathbb{D}(\mu)$ and $\mathbb{D}^{\mu}$ are isomorphic by def\/inition and the f\/irst algebra is studied
in~\cite[Chapter~9]{Dixmier} and~\cite{Lepowsky1973}.
It would be interesting to determine its generators in our situation, where the branching is multiplicity free and the
rank is one.

Since we do not have a~precise description of the algebra $\mathbb{D}(\mu)$, we content to stick to the original
def\/inition of a~MVCP and we provide a~method to f\/ind inf\/initely many examples of them.
To this end we exploit the existence of a~special dif\/ferential operator $\Omega\in\mathbb{D}(\mu)$, the (image of the)
second-order Casimir operator~$\Omega$ on~$G$.
After conjugation with $\Psi_{0}^{\mu}$ we f\/ind a~second-order dif\/ferential operator $D^{\mu}\in\mathbb{D}^{\mu}$ that
is symmetric with respect to the pairing and thus has the MVOPs as simultaneous eigenfunctions.

We show that our method is ef\/fective by applying it to low dimensional examples.
First we classify the data for which we obtain MVCPs of size $2\times2$ from~\cite{HvP}.
The short list that we obtain contains old and new items.
Among the new ones are the symplectic (symmetric) pairs $(\mathrm{USp}(2n),\mathrm{USp}(2n-2)\times\mathrm{USp}(2))$ and
the spherical (but non-symmetric) pair $(\mathrm{G}_{2},\mathrm{SU}(3))$.
For almost all these examples we provide the corresponding MVCPs below.

This paper is organized as follows: In Section~\ref{Section2} we review the construction of families of MVOPs based on the
representation theory of compact Gelfand pairs $(G,K)$.
We restrict the functions $\Psi^{\mu}_{d}$ to the Cartan circle $A\subset G$ and we identify $A=S^{1}$.
The coordinate on~$S^{1}$ is the fundamental zonal spherical function~$\phi$, normalized by two constants~$c$,~$d$ so that
$x=c\phi+d\in[0,1]$.
With this change of variables we denote $\widetilde \Psi_{d}^\mu(x)=\Psi_{d}^\mu(\phi(a))$, which is a~function on~$[0,1]$.

The MVOPs are basically a~ref\/lection on the three term recurrence relations of the func\-tions~$\Psi_d^{\mu}$ and~$\widetilde \Psi_d^{\mu}$, introduced in this section.
More precisely let $\widetilde Q_d$ be def\/ined~by
\begin{gather*}
\widetilde Q_{d}^{\mu}(x)=\big(\widetilde\Psi_0^{\mu}(x)\big)^{-1}\widetilde \Psi_d^{\mu}(x).
\end{gather*}
Then it follows from the three term recurrence relation for $\widetilde\Psi_d^{\mu}$ that $\widetilde Q_{d}^{\mu}$ is
a~polynomial of degree~$d$ with non singular leading coef\/f\/icient.

Section~\ref{Section3} is the heart of this paper.
Here we discuss how the algebra of dif\/ferential opera\-tors~$\mathbb{D}^{\mu}$ comes about.
By means of the bispectral property, we prove that there exist constant matrices $\widetilde R$ and $\widetilde S$ such
that $\widetilde \Psi^{\mu}_0$ satisf\/ies the f\/irst-order dif\/ferential equation
\begin{gather}
\label{eq:intro_first_order_equation_x}
x(1-x)\partial_{x}\widetilde \Psi_0^ \mu(x)=\widetilde \Psi_0^ \mu(x)\big(\widetilde S+x\widetilde R\big).
\end{gather}

\begin{Remark}
If we let $\partial_x$ act on both sides~\eqref{eq:intro_first_order_equation_x}, we obtain an instance of Tirao's
matrix valued dif\/ferential equation~\cite{Tirao-PNAS}.
However, the techniques in~\cite{Tirao-PNAS} do not apply directly because~$\widetilde S$ might have eigenvalues in~$-\mathbb{N}_0$.
This is indeed the case for all the examples considered in this paper.
\end{Remark}

\begin{Remark}
Observe that~\eqref{eq:intro_first_order_equation_x}, can be seen as a~dif\/ferential operator acting on the right on~$\widetilde \Psi^{\mu}_{0}$.
On the other hand~$\widetilde \Psi^{\mu}_{0}$ is also an eigenfunction of the radial part of the Casimir
operator~$\Omega$ of~$G$ acting on the left.
\end{Remark}

In Corollary~\ref{cor:operator_D_tilde}, we exploit again the bispectral property of the functions $\widetilde
\Psi^{\mu}_{d}$ to deduce that the image (radial part) $\widetilde D^{\mu}\in\mathbb{D}^{\mu}$ of the Casimir operator
of~$G$ can be expressed in terms of~$\widetilde R$,~$\widetilde S$ and an additional constant matrix coming from~$(G,K)$.
More precisely we show that the polynomials $\widetilde Q^{\mu}_{d}$ satisfy $\widetilde D^\mu \widetilde
Q^{\mu}_{d}=\widetilde Q^{\mu}_{d} \Lambda_d/(rp^2)$, where
\begin{gather*}
\widetilde D^{\mu}=x(x-1)\partial_{x}^2+\left[\frac{\lambda_1m}{rp^2(M-m)}-2\widetilde S +
x\left(\frac{\lambda_1}{rp^2}-2\widetilde R\right)\right] \partial_{x}+\frac{\Lambda_0}{rp^2},
\end{gather*}
where $M$, $m$ are the maximum and minimum of $\phi|_{S^1}$,~$p$ the period of~$\phi$ and~$r$ a~scaling factor.
This data is provided in Table~\ref{table2} for the various cases.
The diagonal matrix $\Lambda_0$ is the eigenvalue of $\widetilde \Psi^\mu_0$ as an eigenfunction of the Casimir
operator~$\Omega$ and it can be calculated for each pair $(G,K)$.
It follows that the explicit knowledge of the function $\widetilde \Psi_{0}^{\mu}$ implies explicit knowledge of the
corresponding pair $(W^{\mu},\widetilde D^{\mu})$.

In Section~\ref{Section4} we discuss an algorithm to obtain an explicit expression for $\widetilde \Psi_{0}^{\mu}$.
This algorithm can be implemented in \verb|GAP|~\cite{GAP4,vPR} for each specif\/ic pair $(G,K)$.
Once we have a~formula for $\Psi_{0}^{\mu}$ we can calculate the corresponding MVCP by dif\/ferentiation and matrix
multiplication.
We also propose a~method of f\/inding families of MVCPs.
\begin{itemize}\itemsep=0pt
\item[(1)] Take a~family $(G_{n},K_{n},\mu_{n})_{n\in\mathbb{N}}$ for which the construction applies.
For instance, take a~constant family of Gelfand pairs and consider $(G,K,n\mu)$, where $\mu\in F$, or take a~canonical
element of a~face~$F$, for instance the f\/irst fundamental weight~$\omega$, and consider the family
$(G_{n},K_{n},\omega)_{n\in\mathbb{N}}$, where we let the Gelfand pairs $(G_{n},K_{n})_{n\in\mathbb{N}}$ vary with~$n$.
\item[(2)] Calculate the f\/irst so many functions $\widetilde \Psi_{0}^{\mu}$ of the family until a~pattern shows up.
This provides an ansatz for a~family of MVCPs.
\item[(3)] Show that every pair is indeed a~MVCP.
This is not dif\/f\/icult, one needs to check three equations~\cite[Theorem~3.1]{DG}.
It turns out that in many cases the group parameter, which is a~priori discrete, may vary continuously within a~certain range.
\end{itemize}

In Section~\ref{section:SPcase} we discuss the implementation in \verb|GAP|~\cite{GAP4} for the Gelfand pair ($\mathrm{USp}(2n)$,
$\mathrm{USp}(2n-2)\times\mathrm{USp}(2)$)
and the appropriate faces~$F$.
We discuss a~branching rule that is necessary to implement the algorithm.
The branching laws for the symplectic groups are more dif\/f\/icult than those for the special unitary and orthogonal groups.
Indeed, the multiplicities are not only determined by interlacing condition, but also by an alternating sum of partition functions.
At this point it is important to select the right irreducible~$K$-types.
Selecting an irreducible~$K$-representation of highest weight $\mu\in F$, where $(G,K,F)$ is a~multiplicity free system,
guarantees that the branching rules simplify and that the involved algebras are commutative.
In fact, the whole construction of MVOPs would not work for more general irreducible~$K$-representations.

There are two families of MVCPs of size $2\times2$ related to $(\mathrm{USp}(2n),\mathrm{USp}(2n-2)\times\mathrm{USp}(2))$.
We calculate the f\/irst family by hand.
For the other family we calculated the corresponding family of MVCPs using our method.
We apply the machinery once more to give an example of size $3\times3$ that is associated to this Gelfand pair.

In Section~\ref{2x2class} we classify all possible triples $(G,K,\mu)$ that give rise to $2\times2$ MVCPs according to the uniform
construction described in~\cite{HvP}.
Subsequently we determine the corresponding functions $\widetilde \Psi^{\mu}_{0}$.

To indicate that our method is ef\/fective we display most of the MVCPs of size $2\times2$ below.
Some of these MVCPs were already known (Cases~a1, b, d, see~\cite{PT,Pacharoni-Zurrian, TZ}) but they were obtained by dif\/ferent means.
The other MVOPs (Cases~a2, c1, c2, g1, g2, f) are new as far as we know.
The matrix weights are of the form
\begin{gather*}
W^\mu(x)=(1-x)^\alpha x^\beta \widetilde \Psi^{\mu}_{0}(x)^{*}T^{\mu}\widetilde\Psi^{\mu}_{0}(x).
\end{gather*}
We provide the expressions for $\widetilde\Psi^{\mu}_{0}(x)$ and $T^{\mu}$ instead of working out this multiplication,
because the expressions become rather lengthy.
The parameters $n$, $i$, $m$ are a~priori all integers for which we give the bounds in each case, see Remark~\ref{remark:real parameters}.

 {\bf Case a.} The pair $(G,K)=(\mathrm{SU}(n+1),\mathrm{U}(n))$.
We have $\alpha=n-1$, $\beta=0$, $n\ge2$, $1\le i\le n-1$ and $m\in\mathbb{Z}$.
We have two families of MVCPs associated to this example, depending on the sign of~$m$, but in either case
\begin{gather*}
\widetilde D^{\mu}=x(x-1)\partial_{x}^2+\big[{-}1-2\widetilde S + x(n+1-2\widetilde R)\big]
\partial_{x}+\frac{\Lambda_0}{2}.
\end{gather*}

{\bf Case a1.} For $m\geq 0$:
\begin{gather*}
\widetilde \Psi_0^{\mu}(x)=x^{\frac{m}{2}}
\begin{pmatrix}
\sqrt{x} & \sqrt{x}
\vspace{1mm}\\
1 & \dfrac{(m+1)-x(m+n-i+1)}{i-n}
\end{pmatrix},
\qquad
T^{\mu}=
\begin{pmatrix}
i & 0
\\
0 & n-i
\end{pmatrix},
\\
\Lambda_0 =
\begin{pmatrix}
0 & 0
\\
0 & 2(m+n-i+1)
\end{pmatrix},
\qquad
\widetilde R=
\begin{pmatrix}
-\dfrac{m+1}{2} & \dfrac{1}{2}
\vspace{1mm}\\
0 & -\dfrac{m+2}{2}
\end{pmatrix},
\\
\widetilde S=
\begin{pmatrix}
-\dfrac{(m+1)(i-m-n)-m-1-i+n}{2(m+n+1-i)} & \dfrac{m+1}{2(m+n+1-i)}
\vspace{1mm}\\
\dfrac{i-n}{2(i-n-m-1)} & \dfrac{(m+1)(i-m-n)}{2(i-m-n-1)}
\end{pmatrix}.
\end{gather*}

{\bf Case a2.} For $m<0$:
\begin{gather*}
\widetilde \Psi_0^{\mu}(x)=x^{-\frac{(m+1)}{2}}
\begin{pmatrix}
\dfrac{m+(i-m)x}{i} & 1
\vspace{1mm}\\
x^{\frac12} & x^{\frac12}
\end{pmatrix},
\qquad
T^{\mu}=
\begin{pmatrix}
i & 0
\\
0 & n-i
\end{pmatrix},\qquad
\Lambda_0 =
\begin{pmatrix}
0 & 0
\\
0 & m-i
\end{pmatrix},
\\
\widetilde R=
\begin{pmatrix}
\dfrac{m-1}{2} & 0
\vspace{1mm}\\
\dfrac{1}{2} & \dfrac{m}{2}
\end{pmatrix},
\qquad
\widetilde S=
\begin{pmatrix}
-\dfrac{m(i-m-1)}{2(i-m)} & \dfrac{i}{2(i-m)}
\vspace{1mm}\\
-\dfrac{m}{2(i-m)} & -\dfrac{mi+i-m^2}{2(i-m)}
\end{pmatrix}.
\end{gather*}

{\bf Case b.} The pair $(G,K)=(\mathrm{SO}(2n+1),\mathrm{SO}(2n))$.
Let $\alpha=\beta=n-1$, $1\le i\le n-2$.
The corresponding MVCP is given~by
\begin{gather*}
\widetilde D^\mu=x(x-1)\partial_x^2+\big(-2-2\widetilde{S}+2x(n-\widetilde R)\big)\partial_x+\Lambda_0,
\\
\widetilde\Psi_0^{\mu}(x) =
\begin{pmatrix}
2x-1 & 1
\\
1 & 2x-1
\end{pmatrix},
\qquad
T^{\mu}=
\begin{pmatrix}
i & 0
\\
0 & n-i
\end{pmatrix},
\\
\Lambda_{0}=
\begin{pmatrix}
0 & 0
\\
0 & 2(n-i)
\end{pmatrix},
\qquad
\widetilde R =
\begin{pmatrix}
-1 & 0
\\
0 & -1
\end{pmatrix},
\qquad
\widetilde S =
\begin{pmatrix}
\dfrac{1}{2} & \dfrac{1}{2}
\vspace{1mm}\\
\dfrac{1}{2} & \dfrac{1}{2}
\end{pmatrix}.
\end{gather*}

{\bf Case c.} The pair $(G,K)=(\mathrm{USp}(2n),\mathrm{USp}(2n-2)\times\mathrm{USp}(2))$.
We have $\alpha=2n-3$, $\beta=1$ with $n\ge3$.
We have two families of MVCPs associated to this example but in either case
\begin{gather*}
\widetilde D^\mu =x(x-1)\partial_x^2+\big({-}n-2\widetilde{S}+2x(2n-\widetilde R)\big)\partial_x+\Lambda_0.
\end{gather*}

{\bf Case c1:}
\begin{gather*}
\widetilde \Psi_0^{\mu}(x)=
\begin{pmatrix}
\sqrt{x} & \sqrt{x}
\vspace{1mm}\\
1 & \dfrac{x(n-1)-1}{n-2}
\end{pmatrix},
\qquad
T^{\mu}=
\begin{pmatrix}
2 & 0
\\
0 & 2n-4
\end{pmatrix},
\\
\Lambda_0 =
\begin{pmatrix}
0 & 0
\\
0 & 4(n-1)
\end{pmatrix},
\qquad
\widetilde R=
\begin{pmatrix}
-\dfrac12 & \dfrac12
\vspace{1mm}\\
0 & -1
\end{pmatrix},
\qquad
\widetilde S=
\begin{pmatrix}
\dfrac{1}{2(n-1)} & \dfrac{1}{2(n-1)}
\vspace{1mm}\\
\dfrac{(n-2)}{2(n-1)} & \dfrac{(n-2)}{2(n-1)}
\end{pmatrix}.
\end{gather*}

{\bf Case c2:}
\begin{gather*}
\widetilde \Psi^{\mu}_0(x)=
\begin{pmatrix}
\dfrac{x+1}{2} & \dfrac{(n+1)x-2}{n-1}
\vspace{1mm}\\
\sqrt{x} & \dfrac{\sqrt{x}((n+3)x+n-5)}{2(n-1)}
\end{pmatrix},
\qquad
T^{\mu}=
\begin{pmatrix}
\dfrac{2}{n+1} & 0
\vspace{1mm}\\
0 & \dfrac{2}{(n-2)}
\end{pmatrix},
\\
\Lambda_0
\begin{pmatrix}
0 & 0
\\
0 & 2n+6
\end{pmatrix},
\qquad
\widetilde R=
\begin{pmatrix}
-1 & \dfrac{(n+1)}{(n-1)}
\vspace{1mm}\\
0 & -\dfrac{3}{2}
\end{pmatrix},
\qquad
\widetilde S=
\begin{pmatrix}
\dfrac{4}{n+3} & \dfrac{2n-10}{(n-1)(n+3)}
\vspace{1mm}\\
\dfrac{n-1}{n+3} & \dfrac{n-5}{2(n+3)}
\end{pmatrix}.
\end{gather*}

For Case~c2, we do not provide a~formal proof that the family of MVCPs that are produced in this way, are the ones
associated to the Lie theoretical datum.
For this we have to study the various representations, which is quite involved.

 {\bf Case g.} The pair $(G,K)=(\mathrm{G}_{2},\mathrm{SU}(3))$.
There is a~single $2\times 2$ MVCP associated to this pair.
We have $\alpha=\beta=2$ and
\begin{gather*}
\widetilde D^\mu=x(x-1)\partial_x^2+\big({-}3-2\widetilde{S}+x(6-\widetilde R)\big)\partial_x+\frac{\Lambda_0}{2},
\end{gather*}
where
\begin{gather*}
\Psi^\mu(x)=
\begin{pmatrix}
x & x
\\
\sqrt{x} & 3x^{\frac32}-2\sqrt{x}
\end{pmatrix},
\qquad
T^\mu=
\begin{pmatrix}
1 & 0
\\
0 & 2
\end{pmatrix},
\\
\Lambda_0=
\begin{pmatrix}
0 & 0
\\
0 & 6
\end{pmatrix},
\qquad
\widetilde R=
\begin{pmatrix}
-1 & \dfrac{1}{2}
\vspace{1mm}\\
0 & -\dfrac{3}{2}
\end{pmatrix},
\qquad
\widetilde S=
\begin{pmatrix}
\dfrac{5}{6} & \dfrac{1}{3}
\vspace{1mm}\\
\dfrac{1}{6} & \dfrac{2}{3}
\end{pmatrix}.
\end{gather*}

We omit Case~d, $(\mathrm{SO}(2n),\mathrm{SO}(2n-1))$, as it is similar to Case~b.
We make a~few remarks concerning these examples.

\begin{Remark}
\label{remark:real parameters}
The parameters $n$, $m$, $i$ in the various examples may vary in $\mathbb{R}$ rather than in~$\mathbb{N}$, within certain
bounds.
The bounds are determined by the question whether the matrix weight is positive.
To see that the pairs $(W,D)$ remain MVCPs, one has to check the symmetry relations~\eqref{eq:symmetry_equations}.
These expressions are meromorphic in the parameters, so they remain valid.
\end{Remark}

\begin{Remark}
In each case the determinant of the weight is a~product of powers of~$x$ and $(1-x)$ times a~constant.
On the one hand this is quite remarkable, for the weight matrices do not seem to have much structure.
However, it turns out that all the weights that we construct have this property.
This follows from Corollary~\ref{cor:det}, which also settles an earlier Conjecture~1.5.3 of~\cite{vP}.
\end{Remark}

\begin{Remark}
The matrix weights~$W$ may have symmetries, i.e.~they may be conjugated by a~constant matrix into a~diagonal matrix weight.
We check whether this occurs by looking at the commutant of $W^\mu$.
It turns out that we only have non trivial commutant in Cases~b1 and~ b2 for specif\/ic parameters.
\end{Remark}

\section{Lie theoretical background}\label{Section2}

Let $(G,K)$ be a~pair of compact connected Lie groups from Table~\ref{table mfs compact} and let
$\mathfrak{g}$, $\mathfrak{k}$ denote their Lie algebras.
The quotient $G/K$ is a~two-point-homogeneous space (cf.~\cite{Wang}) which implies that~$K$ acts transitively on the
unit sphere in $T_{eK}G/K$.
We denote $T_{eK}G/K=\mathfrak{p}$ and f\/ix a~one dimensional abelian subspace $\mathfrak{a}\subset\mathfrak{p}$.
The one dimensional subspace $\mathfrak{a}\subset\mathfrak{g}$ is the Lie algebra of a~subtorus $A\subset G$ and it
follows that we have a~decomposition $G=KAK$.

Let $M=Z_K(A)$ denote the centralizer of~$A$ in~$K$ with Lie algebra $\mathfrak{m}\subset\mathfrak{k}$.
Let $T_{M}\subset M$ be a~maximal torus and let $T_K\subset K$ be a~maximal torus that extends $T_{M}$.
Then $M\cap T_K=T_M$ and $AT_M$ is a~maximal torus of~$G$.
The Lie algebras of $T_{K}$ and $T_M$ are denoted by $\mathfrak{t}_{K}$ and $\mathfrak{t}_M$.
The complexif\/ications of the Lie algebras are denoted by $\mathfrak{g}_c,\ldots$.

The (restricted) roots of the pairs $(\mathfrak{g}_c,\mathfrak{a}_c\oplus\mathfrak{t}_{M,c})$,
$(\mathfrak{g}_c,\mathfrak{a}_c)$, $(\mathfrak{m}_c,\mathfrak{t}_{M,c})$ and $(\mathfrak{k}_c,\mathfrak{t}_{K,c})$ are
denoted by $R_G$, $R_{(G,A)}$, $R_M$ and $R_K$ respectively.
We f\/ix systems of positive roots $R^+_G$, $R^+_{(G,A)}$, $R^+_M$ and~$R^+_K$ such that the natural projections $R_G\to
R_{(G,A)}$ and $R_G\to R_M$ respect positivity.

The lattices of integral weights of $G$, $K$ and~$M$ are denoted by $P_G$, $P_K$ and $P_M$, the cones of positive integral
weights by $P^+_G$, $P^+_K$ and $P^+_M$.
The theorem of the highest weight implies that the equivalence classes of the irreducible representations are
parametrized by the cones of positive integral weights.
Given $\lambda\in P^+_G$ we denote by $\pi^G_\lambda:G\to\mathrm{GL}(V^{G}_\lambda)$ an irreducible representation of
highest weight~$\lambda$.
The restriction $\pi^G_\lambda|_K$ decomposes into a~f\/inite sum of irreducible~$K$-representations and
for $\mu\in P^+_K$ we denote the multiplicity by $m_\lambda^{G,K}(\mu)=[\pi^G_\lambda|_K:\pi^K_\mu]$.

\begin{Definition}\quad
\begin{itemize}\itemsep=0pt
\item Let $\mu\in P^+_K$.
A~triple $(G,K,\mu)$ is called a~multiplicity free triple if $m_\lambda^{G,K}(\mu)\le1$ for all $\lambda\in P^+_G$.
\item Let $F\subset P^+_K$ be a~face, i.e.~the $\mathbb{N}$-span of some fundamental weights of~$K$.
A~triple $(G,K,F)$ is called a~multiplicity free system if $(G,K,\mu)$ is a~multiplicity free triple for all $\mu\in F$.
\end{itemize}
\end{Definition}

The notion of a~multiplicity free system can be considered for any compact Lie group~$G$ with closed subgroup~$K$.
In~\cite{HvP, vP} it is shown that for $(G,K,F)$ to be a~multiplicity free triple, the pair $(G,K)$ is necessarily
a~Gelfand pair.
Furthermore, the multiplicity free systems $(G,K,F)$ with $(G,K)$ a~Gelfand pair of rank one are classif\/ied by the rows
of Table~\ref{table mfs compact}.
The multiplicity free systems with $(G,K)$ a~compact symmetric pair have been classif\/ied in~\cite{He}.

\begin{table}[ht]
\centering
 \caption{Compact multiplicity free systems of rank one.
In the third column we have given the highest weight $\lambda_{\mathrm{sph}}\in P_{G}^{+}$ of the fundamental zonal
spherical representation in the notation for root systems of Knapp~\cite[Appendix~C]{Knapp}, except for the case
$(G,K)=(\mathrm{SO}_{4}(\mathbb{C}),\mathrm{SO}_{3}(\mathbb{C}))$, where~$G$ is not simple and
$\lambda_{\mathrm{sph}}=\varpi_{1}+\varpi_{2}\in P^{+}_{G}=\mathbb{N}\varpi_{1}+\mathbb{N}\varpi_{2}$.
The groups~$M$ are isogenous to
$\mathrm{U}(n-2)$, $\mathrm{SO}(2n-2)$, $SO(2n-1)$, $\mathrm{USp}(2n-4)\times\mathrm{USp}(2)$, $\mathrm{Spin}(7)$, $\mathrm{SU}(3)$ and
$\mathrm{SU}(2)$ respectively.}\label{table mfs compact}

\vspace{1mm}

\begin{tabular}{|cc|c|c|c|}
\hline
 $G$ && $K$ & $\lambda_{\mathrm{sph}}$ & faces~$F$
\\
\hline
$\mathrm{SU}(n+1)$ & $n\ge1$ & $\mathrm{U}(n)$ & $\varpi_{1}+\varpi_{n}$ & any
\\
$\mathrm{SO}(2n)$ & $n\ge2$ & $\mathrm{SO}(2n-1)$ & $\varpi_{1}$& any
\\
$\mathrm{SO}(2n+1)$ & $n\ge2$ & $\mathrm{SO}(2n)$ & $\varpi_{1}$ & any
\\
$\mathrm{USp}(2n)$ & $n\ge3$ & $\mathrm{USp}(2n-2)\times\mathrm{USp}(2)$ & $\varpi_{2}$& $\dim F\le 2$
\\
\hline
$\mathrm{F}_{4}$ && $\mathrm{Spin}(9)$ & $\varpi_{1}$ & $\dim F\le 1$ or
\\
&&&&$F=\mathbb{N}\omega_{1}+\mathbb{N}\omega_{2}$
\\
\hline
$\mathrm{Spin}(7)$ && $\mathrm{G}_{2}$ & $\varpi_{3}$ & $\dim F\le1$
\\
$\mathrm{G}_{2}$ && $\mathrm{SU}(3)$ & $\varpi_{1}$ & $\dim F\le1$
\\
\hline
\end{tabular}
\end{table}

Let $(G,K,F)$ be a~multiplicity free system from Table~\ref{table mfs compact}, let $\mu\in F$ and def\/ine
$P^+_G(\mu)=\{\lambda\in P^+_G:m_\lambda^{G,K}(\mu)=1\}$.
Let $P^+_M(\mu)=\{\nu\in P^+_M:m_\mu^{K,M}(\nu)=1\}$.
The structure of the set $P^+_G(\mu)$ is important for the construction of the families of matrix valued orthogonal
polynomials that we associate to $(G,K,\mu)$.
In the special case $\mu=0$ we have $P^+_G(0)=\mathbb{N}\lambda_{\mathrm{sph}}$, where $\lambda_{\mathrm{sph}}$ is
called the fundamental spherical weight.
We have indicated the weights $\lambda_{\mathrm{sph}}$ in Table~\ref{table mfs compact}.

The complexif\/ied Lie algebra $\mathfrak{g}_c$ has a~decomposition
$\mathfrak{g}_c=\mathfrak{k}_c\oplus\mathfrak{a}_c\oplus\mathfrak{n}^+$, where $\mathfrak{n}^+$ is the sum of the root
spaces of the positive restricted roots $R_{(G,A)}^+$.
If $(G,K)$ is symmetric this is just the Iwasawa decomposition.
For the two non-symmetric cases see, e.g.,~\cite{vP}.

For $\lambda\in P^+_G(\mu)$ the action of~$M$ on $(V^{G}_\lambda)^{\mathfrak{n}^+}=\{v\in
V^{G}_\lambda:\mathfrak{n}^+v=0\}$ is irreducible of highest weight $\lambda|_{\mathfrak{t}_M}$.
Moreover, $\lambda|_{\mathfrak{t}_M}\in P^+_M(\mu)$, see e.g.~\cite{HvP}.
It follows that the natural projection $q:P_G\to P_M$ induces a~map $P^+_G(\mu)\to P^+_M(\mu)$ which turns out to be
surjective, see~\cite{HvP}.
On the other hand, if $\lambda\in P^+_G(\mu)$ then $\lambda+\lambda_{\mathrm{sph}}\in P^+_G(\mu)$, which follows from an
application of the Borel--Weil theorem.
Def\/ine the degree $d:P^{+}_{G}(\mu)\to\mathbb{N}$~by
\begin{gather*}
d(\lambda+\lambda_{\mathrm{sph}})=d(\lambda)+1,
\qquad
\min\{d(P^{+}_{G}(\mu)\cap(\lambda+\mathbb{Z}\lambda_{\mathrm{sph}}))\}=0.
\end{gather*}
Let $B(\mu)=\{\lambda\in P^+_G(\mu):d(\lambda)=0\}$.
The set $P^+_G(\mu)$ is called the~$\mu$-well and $B(\mu)$ the bottom of the~$\mu$-well.

\begin{Theorem}
\label{thm:parametrization mu-well}
We have $P^+_G(\mu)=\mathbb{N}\lambda_{\mathrm{sph}}+B(\mu)$.
There is an isomorphism
\begin{gather*}
\lambda: \ \mathbb{N}\times P^{+}_{M}(\mu)\to P_{G}^{+}(\mu)
\end{gather*}
such that $\lambda(d,\nu)|_{\mathfrak{t}_M}=\nu$ and $d(\lambda(d,\nu))=d$.
If $\lambda,\lambda'\in P^+_G(\mu)$ and $[\pi^G_{\lambda}\otimes\pi^G_{\lambda_{\mathrm{sph}}}:\pi^G_{\lambda'}]\ge1$
then $d(\lambda)-1\le d(\lambda')\le d(\lambda)+1$.
\end{Theorem}

The proof of Theorem~\ref{thm:parametrization mu-well} is based on a~case by case inspection, see~\cite{HvP,vP}.
Let $\preceq_{\mu}$ be the partial ordering on $P^+_G(\mu)$ induced from the partial ordering on
$\mathbb{N}\lambda_{\mathrm{sph}}+B(\mu)$ given by the lexicographic ordering $(\le,\preceq)$, where $\le$ comes f\/irst.

\begin{Corollary}
\label{cor:upper triangular}
If $\lambda,\lambda'\in P^+_G(\mu)$ and $\big[\pi^G_{\lambda}\otimes\pi^G_{\lambda_{\mathrm{sph}}}:\pi^G_{\lambda'}\big]\ge1$,
then $\lambda-\lambda_{\mathrm{sph}}\preceq_{\mu}\lambda'\preceq_{\mu}\lambda+\lambda_{\mathrm{sph}}$ $($whenever
$\lambda-\lambda_{\mathrm{sph}}\in P^{+}_{G})$.
\end{Corollary}

Indeed, the highest weights of the irreducible representations that occur in the tensor product decompositions are of
the form $\lambda'=\lambda+\lambda''$, where $\lambda''$ is a~weight of $\lambda_{\mathrm{sph}}$ (see, e.g.,~\cite[Proposition~9.72]{Knapp}).
The lowest weight of the fundamental spherical representation is $w_{0}(\lambda_{\mathrm{sph}})$, which equals
$-\lambda_{\mathrm{sph}}$ by inspection.
Here, $w_{0}$ denotes the longest Weyl group element.

Fix a~multiplicity free system $(G,K,F)$ from Table~\ref{table mfs compact} and f\/ix $\mu\in F$.
Let $\pi^K_{\mu}:K\to\mathrm{GL}(V^{K}_{\mu})$ be an irreducible representation of highest weight~$\mu$.
Let $R(G)$ denote the (convolution) algebra of matrix coef\/f\/icients of~$G$ and def\/ine the $(K\times K)$-action on
$R(G)\otimes\End(V^{K}_{\mu})$~by
\begin{gather*}
(k_{1},k_{2})(m\otimes Y)(g)=m\big(k_{1}^{-1}gk_{2}\big)\otimes\pi^{K}_{\mu}(k_{1})Y\pi^{K}_{\mu}(k_{2})^{-1}.
\end{gather*}
The space $E^{\mu}:=(R(G)\otimes\End(V^{K}_{\mu}))^{K\times K}$ is called the space of~$\mu$-spherical
functions.
Note that $\Phi\in E^{\mu}$ satisf\/ies
\begin{gather*}
\Phi(k_{1}gk_{2})=\pi^{K}_{\mu}(k_{1})\Phi(g)\pi^{K}_{\mu}(k_{2})
\qquad
\forall\, k_{1},k_{2}\in K,g\in G.
\end{gather*}
Furthermore, note that $E^0$ is a~polynomial algebra and that $E^{\mu}$ is a~free, f\/initely generated $E^{0}$-module.
In fact, $E^{\mu}\cong E^0\otimes\mathbb{C}^{|B(\mu)|}$ as $E^0$-modules.

Let $\lambda\in P_{G}^{+}(\mu)$ and let $\pi^G_{\lambda}:G\to\mathrm{GL}(V^{G}_{\lambda})$ denote the corresponding
representation.
Then $V_{\lambda}^{G}=V_{\mu}^{K}\oplus (V_{\mu}^{K})^{\perp}$ and we denote by $b:V^{K}_{\mu}\to V^{G}_{\lambda}$
a~unitary~$K$-equivariant embedding and by $b^{*}:V^{G}_{\lambda}\to V^{K}_{\mu}$ its Hermitian adjoint.

\begin{Definition}
The elementary spherical function of type~$\mu$ associated to $\lambda\in P_{G}^{+}(\mu)$ is def\/ined~by
\begin{gather*}
\Phi^{\mu}_{\lambda}: \ G\to\End\big(V^{K}_{\mu}\big):g\mapsto b^{*}\circ\pi_{\lambda}(g)\circ b.
\end{gather*}
\end{Definition}

It is clear that the elementary spherical functions have the desired transformation behaviour.
We equip the space $E^{\mu}$ with a~sesqui-linear form that is linear in the second variable,
\begin{gather*}
\langle\Phi_{1},\Phi_{2}\rangle_{\mu,G}=\int_{G}\tr\left(\Phi_{1}(g)^{*}\Phi_{2}(g)\right)dg
\end{gather*}
with $dg$ the normalized Haar measure.
As a~consequence of Schur orthogonality and the Peter--Weyl theorem we have the following result.

\begin{Theorem}\quad
\begin{itemize}\itemsep=0pt
\item The pairing $\langle\cdot,\cdot\rangle_{\mu,G}:E^{\mu}\times E^{\mu}\to\mathbb{C}$ is a~Hermitian inner product
and $\langle\Phi_{\lambda}^{\mu},\Phi_{\lambda'}^{\mu}\rangle_{\mu,G}=c_{\lambda}\delta_{\lambda,\lambda'}$.
\item $\{\Phi^{\mu}_{\lambda}:\lambda\in P^{+}_{G}(\mu)\}$ is an orthogonal basis of $E^{\mu}$.
\end{itemize}
\end{Theorem}

Denote $\phi=\Phi^{0}_{\lambda_{\mathrm{sph}}}$, the fundamental zonal spherical function and write
\begin{gather}
\label{eq:def_phid}
\phi_{d}=\Phi^{0}_{d\lambda_{\mathrm{sph}}}.
\end{gather}
Then $E^{0}=\mathbb{C}[\phi]$, i.e.~$E^{0}$ is a~polynomial ring generated by~$\phi$.
Note that $\phi(k_{1}gk_{2})=\phi(g)$ for all $k_{1},k_{2}\in K$ and all $g\in G$.
This implies that $\phi\Phi_{\lambda}$ can be expressed as a~linear combination of elementary spherical functions.
It follows from Corollary~\ref{cor:upper triangular} that
\begin{gather}
\label{eqn:recurrence el sf}
\phi\Phi_{\lambda}
=\sum\limits_{\lambda-\lambda_{\mathrm{sph}}\preceq_{\mu}\lambda'\preceq_{\mu}\lambda+\lambda_{\mathrm{sph}}}c_{\lambda,\lambda'}^{\mu}\Phi_{\lambda'}.
\end{gather}
The Borel--Weil theorem implies that $c^{\mu}_{\lambda,\lambda+\lambda_{\mathrm{sph}}}\ne0$ and we can express the
elementary spherical function $\Phi_{\lambda}$ as a~$E^0$-linear combination of the functions $\Phi_{\lambda(0,\nu)}$,
with $\nu\in P_{M}^+(\mu)$.

\begin{Definition}\quad
\begin{itemize}\itemsep=0pt
\item For $\lambda\in P_{G}^{+}(\mu)$ def\/ine $Q_{\lambda}(\phi)=(q^{\mu}_{\lambda,\nu}(\phi):\nu\in B(\mu))$ in
$\mathbb{C}^{|B(\mu)|}[\phi]$~by
\begin{gather*}
\Phi_{\lambda}=\sum\limits_{\nu\in B(\mu)}q^{\mu}_{\lambda,\nu}(\phi)\Phi_{\nu}.
\end{gather*}
\item For $d\in\mathbb{N}$ def\/ine $Q_{d}(\phi)\in\End(\mathbb{C}^{|B(\mu)|})[\phi]$ as the matrix valued
polynomial having the $Q_{\lambda(d,\nu)}(\phi)$ as columns ($\nu\in B(\mu)$).
\end{itemize}
\end{Definition}

\begin{Theorem}
The matrix valued polynomial $Q_{d}$ is of degree~$d$ and has invertible leading coefficient.
\end{Theorem}

Indeed, from~\eqref{eqn:recurrence el sf} we deduce that for each $d\in\mathbb{N}$ there are $A_d$, $B_d$, $C_d$ in
$\End(\mathbb{C}^{|B(\mu)|})$ such that
\begin{gather}
\label{eqn:rec coefs}
\phi Q_d(\phi)=Q_{d+1}(\phi)A_{d}+Q_{d}(\phi)B_{d}+Q_{d-1}(\phi)C_{d}.
\end{gather}
Corollary~\ref{cor:upper triangular} implies that the matrices $A_{d}$ are upper triangular and the non vanishing of
$c^{\mu}_{\lambda,\lambda+\lambda_{\mathrm{sph}}}$ that the diagonals are non-zero.

Def\/ine $V^{\mu}:G\to\End(\mathbb{C}^{|B(\mu)|})$~by
$V^{\mu}(g)_{\nu,\nu'}=\tr(\Phi_{\lambda(0,\nu)}(g)^*\Phi_{\lambda(0,\nu')}(g))$.
We see that $V^{\mu}$ is~$K$-biinvariant, hence it is of the form
$V^{\mu}=\widetilde{W}^{\mu}(\phi)\in\End(\mathbb{C}^{|B(\mu)|})[\phi]$.

The pairing $\langle Q,Q'\rangle=\int_G Q(\phi(g))^*\widetilde{W}^{\mu}(\phi)Q'(\phi(g))dg$ is a~matrix valued inner
product, see e.g.~\cite{HvP}.
Note that all the functions in the integrand are polynomials in~$\phi$, and~$\phi$ is~$K$-biinvariant.
In view of $G=KAK$ and the integral formulas for this decomposition, we have
\begin{gather}
\label{eqn:MVOP pairing}
\langle Q,Q'\rangle=\int_{0}^1 Q(x)^*\widetilde{W}^{\mu}(x)Q'(x)(1-x)^\alpha x^\beta dx,
\end{gather}
where $x=c\phi+d$ (with constants $c$, $d$ depending on the pair $(G,K)$) is a~normalization such that~$x$ attains all
values in $[0,1]$.
The factor $(1-x)^\alpha x^\beta$ is the ordinary Jacobi weight that is associated to the Riemann symmetric space $G/K$
on the interval $[0,1]$.
We denote $W^{\mu}(x)=(1-x)^\alpha x^\beta \widetilde{W}^{\mu}(x)$

Now we come to a~dif\/ferent discription of the family $(Q_d:d\in\mathbb{N})$, one that allows us to transfer
dif\/ferentiability properties of the elementary spherical functions to similar properties of the matrix valued
polynomials.

The spherical functions $\Phi\in E^{\mu}$ are determined by their restriction to~$A$ and
$\Phi(a)\in\End_M(V_{\mu})$ because~$A$ and~$M$ commute.
Since $m^{K,M}_{\mu}(\nu)\le1$, $\End_M(V_{\mu})$ consists of diagonal matrices.
More precisely, with respect to a~basis of the~$M$-subrepresentations of $V_{\mu}$, the matrix $\Phi(a)$ is block
diagonal, and every block is a~constant times the identity matrix of size $\dim\nu$.
Sending such a~matrix to a~vector containing these constant provides an isomorphism
$u:\End_M(V_{\mu})\cong\mathbb{C}^{|B(\mu)|}$.

\begin{Definition}
\label{def:Psi^mu_d}
The function $\Psi^{\mu}_{\lambda}:A\to\mathbb{C}^{|B(\mu)|}$ is def\/ined~by
$\Psi^{\mu}_{\lambda}(a)=u(\Phi^{\mu}_{\lambda}(a))$.
The function $\Psi_d^{\mu}:A\to\End(\mathbb{C}^{|B(\mu)|})$ is the matrix valued function whose columns are the
vector valued functions $\Psi_{\lambda(d,\nu)}^{\mu}$, with $\nu\in P^+_M(\mu)$, and where $\lambda:\mathbb{N}\times
P^{+}_{M}(\mu)\to P^{+}_{G}(\mu)$ is def\/ined in Theorem~\ref{thm:parametrization mu-well}.
\end{Definition}

It follows that $\Psi_d^{\mu}(a)=\Psi_0^{\mu}(a)Q_d(\phi(a))$ and
$\widetilde{W}^{\mu}(\phi(a))=\Psi_0^{\mu}(a)^*T^{\mu}\Psi_0^{\mu}(a)$, with $T^{\mu}$ the diagonal matrix whose entries
are $\dim\nu$.
Moreover, we know precisely which matrix coef\/f\/icients occur in the entries of the functions $\Psi_d^{\mu}$.

\begin{Theorem}
\label{thm:parametrization entries Psi}
The entries of $\Psi_d^{\mu}$ are indexed by the set $P^{+}_{M}(\mu)$.
Hence, the entry $(\Psi_{d}^{\mu})_{\nu_{1},\nu_{2}}$ is equal to the matrix coefficient $m^{\lambda}_{v,v}$, where
$\lambda=\lambda(d,\nu_{2})$ and $v\in V^{M}_{\nu_{1}}\subset V^{K}_{\mu}\subset V^{G}_{\lambda}$ is any vector of
length one.
\end{Theorem}

The proof is immediate from the def\/inition of $\Psi_d^{\mu}$.
Note that the construction of the func\-tions~$\Psi_{d}^{\mu}$ can also be performed for the complexif\/ied pair
$(G_{\mathbb{C}},K_{\mathbb{C}})$.
We obtain a~function on $A_{\mathbb{C}}\cong\mathbb{C}^{\times}$ that takes values in
$\End(\mathbb{C}^{|B(\mu)|})$ and we denote it with the same symbol,
$\Psi_{d}^{\mu}:\mathbb{C}^{\times}\to\End(\mathbb{C}^{N})$.
For later reference, we state the following result concerning the entries of $\Psi_d^{\mu}$, of which the proof is
straightforward.

\begin{Proposition}
\label{prop: degrees}
The entries of the function $\Psi_{\lambda}^{\mu}:A_{\mathbb{C}}\to\mathbb{C}^{|B(\mu)|}$ are Laurent polynomials in
$\mathbb{C}[z]$.
The maximal degree is less than or equal to $|\lambda(H_A)|$, where $H_A$ is defined by $A=\mathfrak{a}/2\pi
iH_A\mathbb{Z}$.
In fact, the maximal degree occurs precisely in the entry labeled with $\lambda|_{\mathfrak{t}_M}$.
\end{Proposition}

The functions $\Psi_{d}^{\mu}$ are analytic, as the entries are matrix coef\/f\/icients.
Moreover, they satisfy the three term recurrence relation
\begin{gather}
\label{eqn:three term}
\phi(z)\Psi_{d}^{\mu}(z)=\Psi_{d+1}^{\mu}(z)A^{\mu}_{d}+\Psi_{d}^{\mu}(z)B^{\mu}_{d}+\Psi_{d-1}^{\mu}(z)C^{\mu}_{d}
\qquad
\text{for all}
\quad
z\in A_{\mathbb{C}},
\end{gather}
where the matrices $A_{d}$, $B_{d}$, $C_{d}$ are the same as in~\eqref{eqn:rec coefs}.
Def\/ine $\Delta^{\mu}(\Psi_{d}^{\mu})=\Psi_{d+1}^{\mu}A_{d}+\Psi_{d}^{\mu}B_{d}+\Psi_{d-1}^{\mu}C_{d}$.
The operator $\Delta^{\mu}$ is a~second-order dif\/ference operator acting on the variable~$d$ that has $\Psi_{d}^{\mu}$
as eigenfunction and with eigenvalue~$\phi$.

\section{Dif\/ferential properties}\label{Section3}

In this section we discuss the second-order dif\/ferential operator that we obtain from the quadratic Casimir
operator~$\Omega$ of the group~$G$.
Moreover, using the fact that $\Psi^{\mu}_{0}$ is an eigenfunction of~$\Omega$ whose eigenvalue is a~diagonal matrix and
that the functions $\Psi^{\mu}_{d}$ satisfy a~three term recurrence relation, we deduce that $\Psi^{\mu}_{0}$ satisf\/ies
a~f\/irst-order dif\/ferential equation.
From the singularities of this equation we deduce that $\Psi^{\mu}_{0}(a)$ is invertible whenever $a\in A$ is a~regular
point for $\phi|_A$.
In this section we work with spherical functions on the complexif\/ied Lie groups $G_{\mathbb{C}}$ and $A_{\mathbb{C}}$.

Let $U(\mathfrak{g}_{c})$ be the universal enveloping algebra of $\mathfrak{g}_{c}$ and let $U(\mathfrak{g}_{c})^{K}$
denote the algebra of $\mathrm{Ad}(K)$-invariant elements.
Let $I(\mu)\subset U(\mathfrak{k}_{c})$ be the kernel of the representation
$U(\mathfrak{k}_{c})\to\End(V_{\mu})$ and def\/ine
\begin{gather*}
\mathbb{D}(\mu):=U(\mathfrak{g}_{c})^{K}\big/\left(U(\mathfrak{g}_{c})^{K}\cap U(\mathfrak{g}_{c})I(\mu)\right).
\end{gather*}
The irreducible representations of $\mathbb{D}(\mu)$ correspond to irreducible representations of $\mathfrak{g}_{c}$
whose restriction to $\mathfrak{k}_{c}$ has a~subrepresentation of highest weight~$\mu$, see e.g.~\cite[Th{\'e}or{\`e}me 9.1.12]{Dixmier}.

The dif\/ferential operators $D\in\mathbb{D}(\mu)$ leave the space of~$\mu$-spherical functions invariant.
As the~$\mu$-spherical functions are determined by their values on~$A$, in view of $G=KAK$ and the transformation
behavior of the~$\mu$-spherical functions, every $D\in\mathbb{D}(\mu)$ def\/ines a~dif\/ferential operator $R(\mu,D)$
satisfying $R(\mu,D)(\Phi^{\mu}|_{A})=(D\Phi^{\mu})|_{A}$ for all $\Phi\in E^{\mu}$.
Since the spherical functions are analytic, we obtain, after identifying $A_{\mathbb{C}}=\mathbb{C}^{\times}$, a~map
$R(\mu):\mathbb{D}(\mu)\to\mathbb{C}(z)\otimes\End(\End_{M}(V^{K}_{\mu}))[\partial_{z}]$, and the map
$R(\mu)$ is an algebra homomorphism~\cite{CM1982}.
The dif\/ferential operator $R(\mu,D)$ is called the radial part of~$D$.
We denote the image of $R(\mu)$ by $\mathbb{D}_{R}(\mu)$.

\begin{Theorem}
For every element $D\in\mathbb{D}_{R}(\mu)$ and every $d\in\mathbb{N}$ there is an element
$\Lambda_{d}(D)\in\End(\mathbb{C}^{N})$ such that
\begin{gather*}
D\Psi^{\mu}_{d}=\Psi^{\mu}_{d}\Lambda_{d}(D).
\end{gather*}
Moreover, the map $\Lambda_{d}:D\mapsto\Lambda_{d}(D)$ is a~representation.
The family of representations $\{\Lambda_{d}\}_{d\in\mathbb{N}}$ separates the points of the algebra $\mathbb{D}(\mu)$.
\end{Theorem}

\begin{proof}
The f\/irst part is proved in~\cite{HvP,vP}.
The matrix $\Lambda_{d}(D)$ is diagonal and its entries are $\Lambda_{d}(D)_{\nu,\nu}=\pi^G_{\lambda(d,\nu)}(D)$.
For the second statement, suppose the converse.
Then there are two dif\/ferential operators $D$, $D'$ that have the same eigenvalues and it follows that $D-D'$ acts as zero
on the elementary spherical functions.
This implies that $D-D'=0$.
\end{proof}

It follows that for any $D\in\mathbb{D}_{R}(\mu)$ and $\Delta^{\mu}$ (def\/ined in~\eqref{eqn:three term}), the triple
$(\Delta^{\mu}, D,\{\Psi^{\mu}_{d}:d\in\mathbb{N}\})$ has a~bispectral property, i.e.~the operators $\Delta^{\mu}$
and~$D$ have the members of the family $\{\Psi^{\mu}_{d}:d\in\mathbb{N}\}$ as simultaneous eigenfunctions (albeit in dif\/ferent variables).
For more on the bispectral property, see e.g.~\cite{DuistermaatGrunbaum, Grunbaum-Tirao}.
In the case $\mu=0$ we have $|B(\mu)|=1$ and we denote the eigenvalues by lower case letters, $D\phi=\lambda(D)\phi$.

The Casimir operator on~$G$ is given as follows.
Let $\{X_{i}:i=1,\ldots,\dim\mathfrak{g}_{c}\}$ be a~basis of $\mathfrak{g}_{c}$ and let
$\{\widetilde{X}_{i}:i=1,\ldots,\dim\mathfrak{g}_{c}\}$ be a~dual basis with respect to the Killing form~$\kappa$ on $\mathfrak{g}_{c}$.
Then $\Omega=\sum\limits_{i,j}\kappa(X_{i},X_{j})\widetilde{X}_{i}\widetilde{X}_{j}$.
On the $\mathfrak{g}_{c}$-representation $V_{\lambda}$ the Casimir operator~$\Omega$ acts with the scalar
$\kappa(\lambda,\lambda)+2\kappa(\lambda,\rho_G)$, where $\rho_{G}$ is half the sum of the positive roots of~$G$.
The image of~$\Omega$ in $\mathbb{D}(\mu)$ is denoted by~$\Omega$ or by $\Omega(\mu)$ if we want to indicate on which
function space we let it act.
Then
\begin{gather}
\label{eq:radial part}
R(\mu,\Omega)=r\big((z\partial_{z})^{2}+c(z)z\partial_{z}+F_{\mu}(z)\big),
\end{gather}
where $c(z)$ is a~meromorphic function on $A_{\mathbb{C}}$, $F_{\mu}(z)$ is a~meromorphic function with matrix
coef\/f\/icients and~$r$ is a~constant that depends on the pair~$(G,K)$.
For the particular case $\mu=0$ we have $F_{0}(z)=0$.
For the symmetric pairs these statements follow from~\cite[Proposition~9.1.2.11]{WarnerII}.
For the two non-symmetric pairs see~\cite{HvP} or~\cite[Paragraph~3.4.28]{vP}.

\begin{Lemma}
\label{lemma: mM}
The function~$\phi$ satisfies $(z\partial_{z}\phi)^{2}=p^2(\phi-M)(\phi-m)$, where~$M$ and~$m$ are the maximum and the
minimum of~$\phi$ restricted to the circle $S^{1}\subset A_{\mathbb{C}}$ and $p=\#(K\cap A)$.
\end{Lemma}

\begin{proof}
The fundamental spherical functions are of the form $\phi(z)=a(z^{p}+z^{-p})/2+b$ (see e.g.~\cite[Table 3.3]{vP}).
On the other hand,~$\phi$ is an eigenfunction of $\Omega(0)$ and the statement follows from a~calculation.
\end{proof}

In what follows we only consider the operator $R(\mu,\Omega)$, the radial part of the image of the second-order Casimir
operator in $\mathbb{D}(\mu)$.
The eigenvalues are denoted by $\Lambda_{d}$ for general $\mu\in F$ and $\lambda_{d}$ for $\mu=0$, i.e.~we have
\begin{gather}
\Omega(0)\phi_{d}=\lambda_{d}\phi_{d},
\label{eqn:eigensystem zero} 
\\
\Omega(\mu)\Psi^{\mu}_{d}=\Psi^{\mu}_{d}\Lambda_{d},
\label{eqn:eigensystem mu}
\end{gather}
where $\phi_{d}$ is given by~\eqref{eq:def_phid}.

\begin{Theorem}
\label{theorem: first order DE for Psi_0}
The function $\Psi^{\mu}_{0}$ satisfies the first-order differential equation
\begin{gather}
\label{eqn:linear system Psi_0}
2rz^{2}\partial_{z}\phi(z)(\partial_{z}\Psi^{\mu}_{0})(z)=\Psi^{\mu}_{0}(z)(R\phi(z)+S),
\end{gather}
on $\mathbb{C}^{\times}$,
where $R=A_{0}^{-1}\Lambda_{1}A_{0}-\Lambda_{0}-\lambda_{1}$ and $S=\Lambda_{0}B_{0}-B_{0}A_{0}^{-1}\Lambda_{1}A_{0}$.
\end{Theorem}
\begin{proof}
We need the identities~\eqref{eqn:three term},~\eqref{eqn:eigensystem zero} and~\eqref{eqn:eigensystem mu} for
$d=1,2$ to prove this result.
Let the operator $r(z\partial_{z})^{2}$ act on both sides of the equality
\begin{gather*}
\phi(z)\Psi^{\mu}_{0}(z)=\Psi^{\mu}_{1}(z)A_{0}+\Psi^{\mu}_{0}(z)B_{0}
\end{gather*}
and work out the dif\/ferentiation.
The derivatives of order ${>}1$ can be written in terms of $\phi$, $\Psi^{\mu}_{0}$ and $\Psi^{\mu}_{1}$ and their f\/irst-order derivatives.
We get
\begin{gather*}
r(z\partial_{z})^{2}(\phi(z)\Psi^{\mu}_{0}(z))=2rz^{2}\phi'(z)(\Psi^{\mu}_{0})'(z)
\\
\qquad{}
+\phi(z)\big[\Psi^{\mu}_{0}(z){\Lambda}_{0}-F_{\mu}(z)\Psi^{\mu}_{0}(z)-c(z)z(\Psi^{\mu}_{0})'(z)\big]+
\big[\lambda_{1}\phi(z)-c(z)z\phi'(z)\big]\Psi^{\mu}_{0}(z),
\end{gather*}
and
\begin{gather*}
r(z\partial_{z})^{2}(\phi(z)\Psi^{\mu}_{0}(z)) =r(z\partial_{z})^{2}\big(\Psi^{\mu}_{1}(z)A_{0}+\Psi^{\mu}_{0}(z)B_{0}\big)
\\
\hphantom{r(z\partial_{z})^{2}(\phi(z)\Psi^{\mu}_{0}(z))}{}
=\left(\Psi^{\mu}_{1}(z){\Lambda}_{1}-F_{\mu}(z)\Psi^{\mu}_{1}(z)-c(z)z(\Psi^{\mu}_{1})'(z)\right)A_{0}
\\
\hphantom{r(z\partial_{z})^{2}(\phi(z)\Psi^{\mu}_{0}(z))=}{}
+\left(\Psi^{\mu}_{0}(z){\Lambda}_{0}-F_{\mu}(z)\Psi^{\mu}_{0}(z)-c(z)z(\Psi^{\mu}_{0})'(z)\right)B_{0}.
\end{gather*}
Equating and using the three term relation and its derivative, we f\/ind
\begin{gather*}
2rz^{2}\phi'(z)(\Psi^{\mu}_{0})'(z)
=\Psi^{\mu}_{1}(z){\Lambda}_{1}A_{0}+\Psi^{\mu}_{0}(z){\Lambda}_{0}B_{0}-\phi(z)\Psi^{\mu}_{0}(z)({\Lambda}_{0}+{\lambda}_{1}).
\end{gather*}
Using the three term recurrence relation once more we get
\begin{gather*}
2rz^{2}\phi'(z)(\Psi^{\mu}_{0})'(z)
\\
\quad{}=\big(\phi(z)\Psi^{\mu}_{0}(z)-\Psi^{\mu}_{0}(z)B_{0}\big)A_{0}^{-1}{\Lambda}_{1}A_{0}
+\Psi^{\mu}_{0}(z){\Lambda}_{0}B_{0}-\phi(z)\Psi^{\mu}_{0}(z)({\Lambda}_{0}+{\lambda}_{1})
\\
\quad{}
=\Psi^{\mu}_{0}(z)\big[\big(A_{0}^{-1}{\Lambda}_{1}A_{0}-{\Lambda}_{0}-{\lambda}_{1}\big)\phi(z)
+{\Lambda}_{0}B_{0}-B_{0}A_{0}^{-1}{\Lambda}_{1}A_{0}\big].
\end{gather*}
Plugging in~$R$ and~$S$ yields the desired equation.
\end{proof}

Note that the matrix~$R$ measures the fact whether the matrices $A_{n}$ of the recurrence relations are diagonal or not.
In the symplectic case these recurrence matrices are not diagonal in general, see Sections~\ref{ss: 2x2-calc}
and~\ref{section:second_2x2_example}.
For the orthogonal groups $(\mathrm{SO}(m+1),\mathrm{SO}(m))$ the recurrence matri\-ces~$A_{n}$ are diagonal.
In this case the matrix~$R$ is diagonal (even scalar in many examples).
For $m=3$ this follows from~\cite[Theorem~3.1]{KvPR2} or~\cite[Theorem~9.4]{Pacharoni-Tirao-Zurrian}.
In general this follows from the description of $P^{+}_{G}(\mu)$ and the decomposition of the tensor product of
a~general irreducible representation and the fundamental spherical representation, see~\cite[Chapter~2.4]{vP}.

Let $\phi_{A_{\mathbb{C}}}$ denote the restriction of~$\phi$ to $A_{\mathbb{C}}$.
Let $A_{\mathbb{C},\mathrm{reg}}$ be the set of points where $\phi_{A_{\mathbb{C}}}$ is an immersion, i.e., where
$d\phi_{A_{\mathbb{C}}}$ is injective.
Denote $A_{\mathrm{reg}}=A\cap A_{\mathbb{C},\mathrm{reg}}$.
Finally let $A_{\mathbb{C},\mu-\mathrm{reg}}$ denote the set of points, where $\det\Psi^{\mu}_0(z)\ne0$.
Since the weight function $W^{\mu}$ is polynomial in~$\phi$ and moreover, the highest degree polynomial occurs precisely
once in each column and once in each row, the determinant of $W^{\mu}$ is a~polynomial in~$\phi$ of positive degree.
It follows that $A_{\mathbb{C},\mu-\mathrm{reg}}\subset A_{\mathbb{C}}$ is a~dense open subset. See also~\cite{vP}.

\begin{Corollary}
\label{cor:det}
If $z\in A_{\mathbb{C},\mathrm{reg}}$ then $\Psi^{\mu}_{0}(z)$ is an invertible matrix.
\end{Corollary}

\begin{proof}
The linear system~\eqref{eqn:linear system Psi_0} has singularities precisely in $A_{\mathbb{C}}\backslash
A_{\mathbb{C},\mathrm{reg}}$.
Hence, locally in $A_{\mathbb{C},\mathrm{reg}}$, there exists a~fundamental solution matrix which is holomorphic and
invertible.
By Theorem~\ref{theorem: first order DE for Psi_0} the function $\Psi_0^{\mu}$ has the same properties on the set
$A_{\mathbb{C},\mu-\mathrm{reg}}$, whose intersection with $A_{\mathbb{C},\mathrm{reg}}$ is dense in $A_{\mathbb{C}}$.
The claim follows.
\end{proof}

Corollary~\ref{cor:det} settles a~conjecture on the determinant of the weight matrices, see~\cite[1.5.3]{vP} and~\cite[Theorem 2.3]{KvPR2}.
Namely, the determinant of $V^{\mu}(\phi)$ is a~polynomial in~$\phi$ that is non-zero outside the critical points of $\phi_{A_{\mathbb{C}}}$.
In view of Lemma~\ref{lemma: mM} $\det(V^{\mu}(\phi))$ is a~multiple times $(\phi-M)^{n_{M}}(\phi-m)^{n_{m}}$, for some $n_{M},n_{m}\in\mathbb{N}$.

We can now conjugate the image $\Omega(\mu)\in\mathbb{D}(\mu)$ of the quadratic Casimir operator~$\Omega$ with the
function $\Psi^\mu_{0}$ to obtain a~dif\/ferential operator of order two for the family of matrix valued orthogonal polynomials.
Since $\Omega(\mu)$ has real eigenvalues the resulting operator is symmetric with respect to the pairing~\eqref{eqn:MVOP pairing}.
Indeed, the matrices $\langle Q_{d},Q_{d}\rangle$ are diagonal matrices.

\begin{Theorem}
\label{theorem:conjugation_Dmu_hypergeometric}
The polynomials $Q_d$, $d\in\mathbb{N}$, satisfy $D^\mu Q_d=Q_d\Lambda_d$, where
\begin{gather*}
D^\mu = (\Psi_0^\mu)^{-1} \Omega \Psi_0^\mu=r(z\partial_{z})^{2}+\big[rc(z)+(R\phi(z)+S)/(z\phi'(z))\big](z\partial_{z})+\Lambda_{0}.
\end{gather*}
\end{Theorem}
\begin{proof}
It is a~straightforward computation that
\begin{gather}
\Omega(\mu)\big(\Psi^\mu_{0}(z)Q(z)\big)= r \Psi^\mu_{0}(z) (z\partial_z)^2 Q(z)
+r\big[2z\partial_z\Psi^\mu_0(z)+c(z)\Psi^\mu_0(z)\big](z\partial_{z}Q(z))
\nonumber
\\
\qquad{}
+r\big[(\partial_{z}\Psi^\mu_0)^2+c(z)(z\partial_{z})\Psi^\mu_0(z)+F_\mu(z)\Psi^\mu_0(z)\big]Q(z).
\label{eq:conjugation_Psi0_proof_thm}
\end{gather}
It follows from~\eqref{eq:radial part} and Lemma~\ref{lemma: mM} that the coef\/f\/icient of $Q(z)$
in~\eqref{eq:conjugation_Psi0_proof_thm} is exactly $\Lambda_0$.
On the other hand, the coef\/f\/icient of $(z\partial_{z}Q(z))$ in~\eqref{eq:conjugation_Psi0_proof_thm} is obtained from
Theorem~\ref{theorem: first order DE for Psi_0}.
Now the theorem follows by multiplying with $(\Psi^\mu_0)^{-1}$ from the left on both sides
of~\eqref{eq:conjugation_Psi0_proof_thm}.
\end{proof}

\begin{Corollary}
\label{cor:operator_D_tilde}
Let $\widetilde \Psi_0^\mu$ be the function on the interval $[0,1]$, defined by $\widetilde
\Psi_0^\mu((\phi(z)-m)/(M-m))=\Psi_0^\mu(z)$, where $z\in\{e^{it}:\, 0\leq t < \pi/p \}$.
Then $\widetilde \Psi_0^\mu$ satisfies the $($right$)$-first-order differential equation
\begin{gather}
\label{eq:first_order_equation_x}
x(1-x)\partial_{x}\widetilde \Psi_0^ \mu(x)=\widetilde{\Psi}_0^ \mu(x)\big(\widetilde S+x\widetilde R\big),
\qquad
\widetilde S=-\frac{S+mR}{2rp^2(M-m)},
\qquad
\widetilde R=-\frac{R}{2rp^2}.
\end{gather}
Let $\widetilde Q_d=Q_d\circ ((\phi-m)/(M-m))$.
Then $\widetilde D^\mu\widetilde{Q}_d=\widetilde{Q}_d (\Lambda_d/(rp^2))$, where $\widetilde D^\mu$ is the differential
operator
\begin{gather}
\label{eq:Dmu}
\widetilde D^{\mu}=x(x-1)\partial_{x}^2+\left[\frac{\lambda_1m}{rp^2(M-m)}-2\widetilde S +
x\left(\frac{\lambda_1}{rp^2}-2\widetilde R\right)\right] \partial_{x}+\frac{\Lambda_0}{rp^2}.
\end{gather}
Here $M$, $m$ are the maximum and minimum of $\phi|_{S^1}$, see Lemma~{\rm \ref{lemma: mM}}.
\end{Corollary}
\begin{proof}
The proof is a~straightforward consequence of Lemma~\ref{lemma: mM},
Theorem~\ref{theorem:conjugation_Dmu_hypergeometric} and~\eqref{eqn:eigensystem zero}.
\end{proof}

\begin{Remark}
We observe that the function $(\phi-m)/(M-m)$ is a~bijection from $\{e^{it}:\, 0\leq t < \pi/p \}$ onto the interval
$[0,1]$ (see Table~\ref{table2}), so that $\widetilde \Psi_0^\mu$ in Corollary~\ref{cor:operator_D_tilde} is well
def\/ined.
\end{Remark}

In Section~\ref{subsection:calculation data} we provide all the data $r$, $p$.
The scaling~$r$ is determined by comparing the radial part of $\Omega(0)$ to the hypergeometric dif\/ferential operator
for the Jacobi polynomials on~$G/K$.
Basically we only need to know $(\alpha,\beta)$.

\begin{Remark}
The operator $\widetilde D^\mu$ in~\eqref{eq:Dmu} is a~matrix valued hypergeometric equation~\cite{Tirao-PNAS}.
In the scalar case, the polynomial eigenfunctions of~\eqref{eq:Dmu} can be written in terms of hypergeometric series.
In the matrix valued setting, in order to give a~simple expression of $\widetilde Q_d$ as matrix valued hypergeometric
series it is necessary to perform a~deeper analysis.
See Corollary~\ref{cor:sf_as_Hypergeometric} for the case of the symplectic group.
\end{Remark}

\begin{Remark}
Theorem~\ref{theorem:conjugation_Dmu_hypergeometric} allows us to calculate the conjugation of the Casimir operator with
the function $\Psi^{\mu}_0$ for individual cases, without calculating the radial part of the Casimir operator.
The latter is in general very technical so Theorem~\ref{theorem:conjugation_Dmu_hypergeometric} makes it much easier to
generate explicit examples of dif\/ferential operators.
For an explicit expression of the dif\/ferential operator $D^{\mu}$ we only need to know the eigenvalue~$\lambda_{1}$ and
an explicit expression of the function~$\Psi^{\mu}_0$.
Indeed, using Lemma~\ref{theorem: first order DE for Psi_0} we f\/ind expressions for~$R$ and~$S$, using computer algebra,
and these, together with~$\lambda_{1}$ gives an expression for $D^{\mu}$ by~\eqref{eq:Dmu}.

We conclude that a~numerical expression of the functions $\Psi^{\mu}_0$ allows us, using computer algebra, to generate
examples of a~matrix valued classical pairs $(W^{\mu},D^{\mu})$.
\end{Remark}

\begin{Remark}
\label{remark:list_of_operators}
The operator $\widetilde D^\mu$ in Corollary~\ref{cor:operator_D_tilde} is given explicitly by the following expressions
(use Table~\ref{subsection:calculation data}).
\begin{itemize}\itemsep=0pt
\item For $\mathrm{SU}(n+1)$: $\widetilde D^{\mu}=x(x-1)\partial_{x}^2+\big[{-}1-2\widetilde S + x(n+1-2\widetilde
R)\big] \partial_{x}+\frac{\Lambda_0}{2}$.
\item For $\mathrm{SO}(2n)$: $\widetilde D^{\mu}=x(x-1)\partial_{x}^2+\big[{-}\frac{2n-1}{2}-2\widetilde S +
x(2n-1-2\widetilde R)\big] \partial_{x}+\Lambda_0$.
\item For $\mathrm{SO}(2n+1)$: $\widetilde D^{\mu}=x(x-1)\partial_{x}^2+\big[{-}n-2\widetilde S + x(2n-2\widetilde
R)\big] \partial_{x}+\Lambda_0$.
\item For $\mathrm{USp}(2n)$: $\widetilde D^{\mu}=x(x-1)\partial_{x}^2+\big[{-}2-2\widetilde S + x(2n-2\widetilde
R)\big] \partial_{x}+\frac{\Lambda_0}{2}$.
\item For $\mathrm{F}_4$: $\widetilde D^{\mu}=x(x-1)\partial_{x}^2+\big[{-}6-2\widetilde S + x(12-2\widetilde R)\big]
\partial_{x}+\Lambda_0$.
\item For $\mathrm{Spin}(7)$: $\widetilde D^{\mu}=x(x-1)\partial_{x}^2+\big[{-}\frac72-2\widetilde S + x(7-2\widetilde
R)\big] \partial_{x}+\frac43 \Lambda_0$.
\item For $\mathrm{G}_2$: $\widetilde D^{\mu}=x(x-1)\partial_{x}^2+\big[{-}3-2\widetilde S + x(6-2\widetilde R)\big]
\partial_{x}+\frac{\Lambda_0}{2}$.
\end{itemize}
\end{Remark}

\section{A method to calculate MVCPs}\label{Section4}

In this section we describe an algorithm to calculate the functions $\Psi^{\mu}_{d}$ from Def\/inition~\ref{def:Psi^mu_d}.
We have implemented this algorithm in the computer package \verb|GAP| for the symmetric pair
$(\mathrm{USp}(2n),\mathrm{USp}(2n-2)\times\mathrm{USp}(2))$ (see Section~\ref{section:SPcase} for this pair
and~\cite{vPR} for the \verb|GAP|-f\/iles).
From the description of the algorithm it is clear how to implement it for other multiplicity free triples in
Table~\ref{table mfs compact}.
In the Section~\ref{S The algorithm} we discuss the general algorithm to calculate $\Psi^{\mu}_{n}$.
In Section~\ref{S Implementation} discuss the implementation of the algorithm in \verb|GAP|.
Finally we provide the necessary data from the compact Gelfand pairs to calculate the MVCPs.

\subsection{The algorithm}\label{S The algorithm}

The algorithm to calculate the functions $\Psi_{d}^{\mu}$ has the following input: (1) a~multiplicity free system
$(G,K,F)$ of rank one from Table~\ref{table mfs compact}, (2) an element $\mu\in F$ and (3) an integer~$d$.
The output is an expression for the function $\Psi_{d}^{\mu}:S^{1}\to\End(\mathbb{C}^{N_{\mu}})$, where
$N_{\mu}$ is the number of elements in the bottom $B(\mu)$ of the~$\mu$-well.

\begin{center}
\begin{tabular}{llr}
\hline
\multicolumn{2}{c}{{\bf The Algorithm}}
\\
\hline
1.
& Determine $P^+_M(\mu)=\{\nu_1,\ldots,\nu_{N_{\mu}}\}$.
\\
2.
& Determine the bottom $B(\mu)=\{b(\mu,\nu_{1}),\ldots, b(\mu,\nu_{N_{\mu}})\}$ of the~$\mu$-well.
\\
3.
& Determine $B^d(\mu):=\{\lambda_j=b(\mu,\nu_{j})+d\lambda_{\mathrm{sph}}:\; j=1,\ldots,N_{\mu}\}$.
\\
4.
& {\bf for all $\lambda_j\in B^d(\mu)$ do}
\\
5.
& \hspace{0.5cm} Determine a~$K$-equivariant embedding $\gamma:V_\mu\to V_{\lambda_j}$.
\\
6.
& \hspace{0.5cm} {\bf for all $\nu_i\in P^+_M(\mu)$ do}
\\
7.
& \hspace{1cm} Determine a~vector $v_{\lambda_j,\mu,\nu_i}\in \gamma(V_\mu)$ of highest weight $\nu_i$
\\
& \hspace{1cm} of length one.
\\
8.
& \hspace{1cm} Determine the matrix coef\/f\/icient $e^{it}\mapsto\langle
a_{t}v_{\lambda_j,\mu,\nu_i},v_{\lambda_j,\mu,\nu_i}\rangle$.
\\
9.
& \hspace{0.5cm} Put all the entries in the~$j$-th column
\\
10.
& Build $\Psi_d^\mu$ with all the columns $j=1,\ldots,N_{\mu}$.
\\
\hline
\end{tabular}
\end{center}

\textbf{Step 1.} Determine $P^{+}_{M}(\mu)$.
The set $P^{+}_{M}(\mu)$ parametrizes the elementary spherical functions of a~f\/ixed degree~$d$
(Theorem~\ref{thm:parametrization mu-well}) and also the entries of the elementary spherical functions restricted to~$A$
(Theorem~\ref{thm:parametrization entries Psi}).
Thus, the number of elements~$N$ in $P^{+}_{M}(\mu)$ determines the size of the matrices $\Psi_{d}^{\mu}(a)$, $a\in A$.

\textbf{Steps 2 and 3.} Determine $B^d(\mu)$.
This set determines the spectrum $P^{+}_{G}(\mu)$ and in turn parametrizes the elementary spherical functions.

\textbf{Step 5.} Determine a~$K$-equivariant embedding $\gamma:V^{K}_{\mu}\to V^{G}_{\lambda_j}$.
To this end we calculate the weight spaces $V^{G}_{\lambda_j}(\mu')$, where $\mu'$ is a~weight in $P_{G}$ that projects
onto $\mu\in P^{+}_{K}$, i.e.~with $q(\mu')=\mu$.
If $\mathrm{rank} G=\mathrm{rank} K$ then $\mu'=\mu$.
The highest weight vector $v_{0}\in V^{K}_{\mu}$ is annihilated by all the simple root vectors of~$K$ and the image
$\gamma(v_{0})$ is the unique vector with this property.
Hence we calculate
\begin{gather}
\label{eqn:mu highest weight line in lambda}
\bigcap_{\mu': q(\mu')=\mu}V^{G}_{\lambda_j}(\mu')\cap\bigcap_{\beta\in\Pi_{K}}\ker(E_{\beta}),
\end{gather}
which is a~one-dimensional subspace of $V_{\lambda_j}$ because $m^{G,K}_{\lambda_j}(\mu)=1$.
Fix a~non-zero element $v_{0}'$ in~\eqref{eqn:mu highest weight line in lambda}.
Put $\gamma(v_{0})=v_{0}'$ and def\/ine $\gamma:V^{K}_{\mu}\to V^{G}_{\lambda_j}$ by stipulating that~$\gamma$
is~$K$-equivariant.
Letting the root vectors in $\mathfrak{k}_{\mathbb{C}}$ of the negative simple roots act on $v_{0}'$ gives a~basis of
$\gamma(V^{K}_{\mu})\subset V^{G}_{\lambda_j}$.

\textbf{Step 7.} Determine a~vector $v_{\lambda_j,\mu,\nu_i}$ in $\gamma(V^{K}_{\mu})$ of highest~$M$-weight $\nu_{i}$.
This is similar to the f\/irst part of  Step~5.

\textbf{Step 8.} Determine the matrix coef\/f\/icient $e^{it}\mapsto\langle
a_{t}v_{\lambda_j,\mu,\nu_i},v_{\lambda_j,\mu,\nu_i}\rangle$.
By Proposition~\ref{prop: degrees} we know that any matrix coef\/f\/icient of $\pi^{G}_{\lambda_j}$ restricted to~$A$ is
a~Laurent polynomial whose terms are of degree~$k$ with $|k|\le |\lambda(H_{A})|$.
On the other hand we have $a_{t} = \exp(iH_{A}t)$, which implies
\begin{gather*}
\langle a_{t}v_{\lambda_j,\mu,\nu_i},v_{\lambda_j,\mu,\nu_i}\rangle_{V_{\lambda}} =
\sum\limits_{p=0}^{\infty}\frac{(it)^{p}}{p!}\langle
H_{A}^{p}v_{\lambda_j,\mu,\nu_{i}},v_{\lambda_j,\mu,\nu_{i}}\rangle_{V^{G}_{\lambda_j}}.
\end{gather*}
It suf\/f\/ices to calculate $\langle H_{A}^{k}v_{\lambda_j,\mu,\nu_{i}},v_{\lambda_j,\mu,\nu_{i}}\rangle_{V^{G}_{\lambda}}$
for $k\le2|\lambda(H_{A})|$.
Indeed, the Laurent polyno\-mial~$m^{\lambda}_{v,w}|_{A}$ is determined by the the f\/irst $2|\lambda(H_{A})|$ terms of its
power series.
To get an invariant inner product on the representation spaces we have implemented the Shapovalov form, see
e.g.~\cite{Shapovalov}.

\textbf{Steps 9 and 10.} The matrix coef\/f\/icients we obtain in  Step~8 are put in a~column vector of length
$N_{\mu}$ and these $N_{\mu}$ column vectors are in turn put as columns in an $(N_{\mu}\times N_{\mu})$-matrix.

\subsection{Implementation in GAP}\label{S Implementation}

In  Step~1 we need to determine the $P^{+}_{M}(\mu)$.
Descriptions of $P^{+}_{M}(\mu)$ are given by classical branching rules for the lines~1,~2 and 3, see~\cite{vP}.
For line 4 see Section~\ref{section:SPcase}, for line~5 see~\cite{Baldoni-Silva} and for the two remaining lines
see~\cite{HvP}.
Upon choosing suitable tori, the branching rules amount to interlacing conditions of strings of integers which are
easily implemented in \verb|GAP|.

To implement the bottom $B(\mu)$ for Steps 2 and 3 we need a~precise description.
For the f\/irst three lines of Table~\ref{table mfs compact} see~\cite{Camporesi2, vP} and for the other lines
see~\cite{HvP}.
Once $B^{d}(\mu)$ is known for $d=0$, it is known for all~$d$.
The bottom depends (piecewise) af\/f\/ine linearly on $P^{+}_{M}(\mu)$ and is thus easily implemented in \verb|GAP|.

Having the highest weights of the involved irreducible $G$, $K$ and~$M$-representations we can do all the Lie algebra
calculations in the appropriate~$G$-representation space $V^{G}_{\lambda}$ using \verb|GAP|.
This settles Steps~5 and~7.

The linear system that we need to solve in  Step~8 to determine the coef\/f\/icients of the Laurent polynomials
$m^{\lambda}_{v_{\nu_{i}},v_{\nu_{i}}}|_{A}$ is easily implemented.
The actual output of the algorithm is a~number of $N^{2}$ strings of real numbers, each string representing a~Laurent
polynomial $m^{\lambda}_{v_{\nu_{i}},v_{\nu_{i}}}|_{A}$.

A simple loop provides the implementation of Steps~9 and~10.

\begin{Remark}
The representation spaces that are used soon become very large, which makes the implementation only suitable for
relatively small~$K$-types~$\mu$.
In this case, \textit{small} means that the calculation is performed in a~reasonable time.
It depends on, among other things, the dimensions of the representations spaces $V^{G}_{\lambda}$, $\lambda\in B(\mu)$.
It would be interesting to make these things more precise, once the algorithm is improved.
For instance, of the spaces $V^{G}_{\lambda}$, $\lambda\in B(\mu)$, we only need a~few vectors, whereas in our
implementation we f\/irst need the whole space $V^{G}_{\lambda}$ to calculate these few vectors.
Unfortunately, we do not see how we can realize these kind of improvements now.
\end{Remark}

\begin{table}[ht]
\centering
\caption{Data needed to calculate MVCPs.}\label{table2}   \vspace{1mm}

\begin{tabular}{|c|c|c|c|c|c|}
\hline
~$G$ &~$K$ & $\lambda_{1}$ & $\phi(a_{t})$ & $(\alpha,\beta)$ &~$r$
\\
\hline
$\mathrm{SU}(n+1)$ & $\mathrm{U}(n)$ & $2n+2$ & $\frac{(n+1)\cos^{2}(t)-1}{n}$ & $(n-1,0)$ & $2p^{-2}$
\\
$\mathrm{SO}(2n)$ & $\mathrm{SO}(2n-1)$ & $2n-1$ & $\cos(t)$ & $(n-\frac{3}{2},n-\frac{3}{2})$ & $p^{-2}$
\\
$\mathrm{SO}(2n+1)$ & $\mathrm{SO}(2n)$ & $2n$ & $\cos(t)$ &$(n-1,n-1)$ & $p^{-2}$
\\
$\mathrm{USp}(2n)$ & $\mathrm{USp}(2n-2)\times\mathrm{USp}(2)$ & $4n$ & $\frac{n\cos^{2}(t)-1}{n-1}$&$(2n-3,1)$ &
$2p^{-2}$
\\
\hline
$\mathrm{F}_{4}$ & $\mathrm{Spin}(9)$ & 12 & $\cos(2t)$& $(7,3)$& $p^{-2}$
\\
\hline
$\mathrm{Spin}(7)$ & $\mathrm{G}_{2}$ & 21/4 & $\cos(3t)$ & $(\frac{5}{2},\frac{5}{2})$& $\frac{3}{4}p^{-2}$
\\
$\mathrm{G}_{2}$ & $\mathrm{SU}(3)$ & 12 & $\cos(2t)$ &$(2,2)$& $2p^{-2}$
\\
\hline
\end{tabular}
\end{table}

\subsection{Obtaining the MVCPs}\label{subsection:calculation data}

Once we have constructed the function $\Psi^{\mu}_{0}$ whose entries are polynomials in $\cos(t)$ we follow the next
steps
\begin{enumerate}\itemsep=0pt
\item[\rm (a)] We make the change of variables $x=(\cos(pt)+1)/2$ so that $\widetilde
\Psi^{\mu}_{0}(x)=\Psi^{\mu}_{0}(t)$.
\item[\rm (b)] We calculate $x(1-x)\big(\widetilde \Psi_0^\mu(x)\big)^{-1}\partial_{x}\widetilde \Psi_0^ \mu(x)$ which gives the
matrices $\widetilde R$, $\widetilde S$ according to Corollary~\ref{cor:operator_D_tilde}.
\item[\rm (c)] The dif\/ferential operator $\widetilde D^{\mu}$ is now given by Remark~\ref{remark:list_of_operators}.
The missing data $m$, $M$, $p$ is contained in the fundamental spherical function~$\phi$ (see the proof of Lemma~\ref{lemma:
mM}), which, as well as~$\lambda_{1}$, we provide in Table~\ref{table2}.
The eigenvalue~$\Lambda_{0}$ has to be calculated in the individual cases using for example Weyl's dimension formula.
We do this in the Sections~\ref{section:SPcase} and~\ref{2x2class} for the $(2\times2)$-examples.
\item[\rm (d)] To obtain an expression for the weight $W^{\mu}(x)=\widetilde{W}^{\mu}(x)w(x)$ we need the scalar weight
$w(x)=(1-x)^{\alpha}x^{\beta}$ and the matrix $T^{\mu}$ that has the dimensions of the irreducible~$M$-representations
$\nu\in P^{+}_{M}(\mu)$ on the diagonal, see~\eqref{eqn:MVOP pairing} and the discussion following Def\/ini\-tion~\ref{def:Psi^mu_d}.
The parameters $(\alpha,\beta)$ of the scalar weights are provided in Table~\ref{table2}.
The weight is given~by
\begin{gather*}
W^{\mu}(x)=(1-x)^{\alpha}x^{\beta}\,\Psi^{\mu}_{0}(x)^{*}T^{\mu}\Psi^{\mu}_{0}(x),
\qquad
x=(\phi-m)/(M-m),
\end{gather*}
Hence, knowledge of $(\Psi^{\mu}_{0},\phi,T^{\mu},\lambda_{1},\Lambda_{0})$ allows us to f\/ind explicit expressions of
the pair $(W^{\mu}, \widetilde D^{\mu})$.
\item[\rm (e)] A~formal proof that $(W^\mu,\widetilde D^\mu)$ is a~MVCP can be done showing that $\widetilde D^\mu$ is
symmetric with respect to $W^{\mu}$.
This boils down to prove that the following symmetry equations,~\cite[Theorem 3.1]{DG}, hold true
\begin{gather}
\label{eq:symmetry_equations}
A_1^*W  = -WA_1 + 2(x(x-1)W)',
\qquad \!\!
A_0^*W  = WA_0 - (WA_1)' +(x(x-1)W)'',\!\!\!
\end{gather}
where
\begin{gather*}
A_1(x)=\frac{\lambda_1m}{rp^2(M-m)}-2\widetilde S + x\left(\frac{\lambda_1}{rp^2}-2\widetilde R\right),
\qquad
A_0=\frac{\Lambda_0}{rp^2}.
\end{gather*}
We observe that the boundary conditions in~\cite[Theorem 3.1]{DG} are always satisf\/ied in our case.
\end{enumerate}

\section{The symplectic case}\label{section:SPcase}

Let $n\ge3$ and $G=\mathrm{USp}(2n)$, $K=\mathrm{USp}(2n-2)\times\mathrm{USp}(2)$.
The branching rules for~$G$ to~$K$ are due to Lepowsky~\cite{Lepowsky}.
Let $K_{1}=\mathrm{USp}(2)\times\mathrm{USp}(2n-4)\times\mathrm{USp}(2)\subset K$ and
$M=\mathrm{USp}(2)\times\mathrm{USp}(2n-4)\subset K_{1}$, where the embedding $K_{1}\subset K$ is the canonical one and
where $M\subset K_{1}$ is given by $(x,y)\mapsto(x,y,x)$.
The branching rules for~$K$ to~$M$ are due to Baldoni-Silva~\cite{Baldoni-Silva}.
If we choose the~$K$-type in a~2-dimensional face then the branching rules become considerably simple.
We employ the standard choices for roots and weights of the symplectic group~\cite[Appendix~C]{Knapp}.
We have $\mathrm{rank} G=\mathrm{rank} K=\mathrm{rank} K_{1}=n$ and $\mathrm{rank} M=n-1$.
The weight lattices of $G$, $K$, $K_{1}$ are equal to $P=\mathbb{Z}^{n}$.
Let $\{\epsilon_{i}:i=1,\ldots,n\}$ denote the standard basis of~$P$.
The fundamental weights of~$G$ are $\varpi_{i}=\sum\limits_{j=1}^{i}\epsilon_{j}$ for $i=1,\ldots,n$.
The fundamental weights of~$K$ are $\omega_{i}=\sum\limits_{j=1}^{i}\epsilon_{j}$ for $i=1,\ldots,n-1$ and
$\omega_{n}=\epsilon_{n}$.
The fundamental weights of $K_{1}$ are $\xi_{1}=\epsilon_{1}$, $\xi_{i}=\sum\limits_{j=2}^{i}\epsilon_{j}$ for
$i=2,\ldots,n-1$ and $\xi_{n}=\epsilon_{n}$.
The fundamental weights of~$M$ are identif\/ied with $\eta_{1}=\frac{1}{2}(\epsilon_{1}+\epsilon_{n})$ and
$\eta_{i}=\sum\limits_{j=2}^{i}\epsilon_{j}$ for $i=2,\ldots,n-1$.
The fundamental weights generate the semi group of integral dominant weights $P^{+}_{K}$, $P^{+}_{K_{1}}$ and $P^{+}_{M}$
which in turn parametrize the equivalence classes of irreducible representations of $K$, $K_{1}$ and~$M$.

\subsection[Branching rule from~$K$ to~$M$]{Branching rule from~$\boldsymbol{K}$ to~$\boldsymbol{M}$}

The branching rule from~$K$ to~$M$ can be calculated in two steps: f\/irst from~$K$ to~$K_{1}$, then from~$K_{1}$ to~$M$.
A~crucial ingredient of the branching rule is the partition function $p_{\Xi}:\mathbb{Z}^{n-1}\to\mathbb{N}$, where
$\Xi=\{\epsilon_{i}\pm\epsilon_{1}:i=2,\ldots,n-1\}$ and
\begin{gather}
\label{def:partition function}
p_{\Xi}(z)=\#\left\{(n_{\xi})_{\xi}\in\mathbb{N}^{|\Xi|}:\sum\limits_{\xi\in\Xi}n_{\xi}\xi=z\right\}.
\end{gather}
\begin{Proposition}
\label{prop: branching}
Let $\mu=x\omega_{i}+y\omega_{j}$ with $i<j$ and $x,y\in\mathbb{N}$.
Write $\mu=\sum\limits_{k=1}^{n}b_{k}\epsilon_{k}$ and let $\nu=\sum\limits_{k=1}^{n}c_{k}\epsilon_{k}\in
P^{+}_{K_{1}}$.
Define $C_{1}=b_{1}-\max(b_{2},c_{2})$, $C_{k}=\min(b_{k},c_{k})-\max(b_{k+1},c_{k+1})$ for $k=2,\ldots,n-2$ and
$C_{n-1}=\min(b_{n-1},c_{n-1})$ and $C_{n}=0$.
We have $m^{K,K_{1}}_{\mu}(\nu)=1$ if and only if
\begin{itemize}\itemsep=0pt
\item[$(1)$] $c_{n}=b_{n}$,
\item[$(2)$] $C_{i}+C_{j}-c_{1}$ is even,
\item[$(3)$] $C_{i}+C_{j}\ge c_{1}\ge |C_{i}-C_{j}|$,
\item[$(4)$] $C_{k}\ge0$ for $k=1,\ldots,n-2$,
\end{itemize}
\end{Proposition}
\begin{proof}
The statement is Lepowsky's branching rule for $\mathrm{USp}(2n-2)$ to $\mathrm{USp}(2)\times\mathrm{USp}(2n-4)$ with
the additional information that we have control over the multiplicity.
The support of the function $m^{K,K_{1}}_{\mu}:P^{+}_{K_{1}}\to\mathbb{N}$ is contained in the set that is determined~by
condition (1) and (4) and the additional condition that $\sum\limits_{k=1}^{n-1}C_{k}$ is even,
see~\cite[Theorem~2]{Lepowsky}.
In this case the multiplicity is given~by
\begin{gather}
m^{K,K_{1}}_{\mu}(\nu)=p_{\Xi}((C_{1}-c_{1})\epsilon_{1}+C_{2}\epsilon_{2}+\dots+C_{n-1}\epsilon_{n-1})
\nonumber
\\
\phantom{m^{K,K_{1}}_{\mu}(\nu)=}{}
-p_{\Xi}((C_{1}+c_{1}+2)\epsilon_{1}+C_{2}\epsilon_{2}+\dots+C_{n-1}\epsilon_{n-1}),
\label{eqn:multiplicity}
\end{gather}
where $p_{\Xi}$ is the partition function~\eqref{def:partition function}.
See~\cite[Theorem~9.50]{Knapp} for a~proof of an equivalent statement (the roots are permuted because the embedding is
dif\/ferent).
If $m^{K,K_{1}}_{\mu}(\nu)=1$ then~$\nu$ is in the support of $m^{K,K_{1}}_{\mu}$ which implies that $C_{k}\ge0$ for
$k=1,\ldots,n-2$.
The condition on~$\mu$ implies that $C_{k}=0$ unless $k=i$ or $k=j$, hence~(2) is satisf\/ied.
Condition~(1) is trivially satisf\/ied.
It remains to check~(3), which we do below.
Conversely, if (1)--(4) are satisf\/ied, then~$\nu$ is in the support of $m^{K_{1},K}_{\mu}$, because $C_{k}=0$ unless
$k=i$ or $k=j$.
It remains to show that~(3) implies $m^{K,K_{1}}_{\mu}(\mu)=1$, which we check by calculating~\eqref{eqn:multiplicity}.
We distinguish two cases: $i=1$ and $i>1$.
Suppose that $i=1$.
Then $C_{k}=0$ unless $k=1$ or $k=j$.
Without loss of generality we may assume that $j=2$ (if $j=n$ then we take $C_{2}=0$) and
$\Xi=\{\epsilon_{2}\pm\epsilon_{1}\}$.
Then
\begin{gather*}
m^{K,K_{1}}_{\mu}(\nu)=p_{\Xi}((C_{1}-c_{1})\epsilon_{1}+C_{2}\epsilon_{2})-p_{\Xi}((C_{1}+c_{1}+2)\epsilon_{1}+C_{2}\epsilon_{2}).
\end{gather*}
Both terms are zero or one, because $|\Xi|=2$.
The f\/irst term is one if and only if $C_{1}-C_{2}\le c_{1}\le C_{1}+C_{2}$ and $C_{1}+C_{2}-c_{1}$ even.
The second term is minus one if and only if $c_{1}\le C_{2}-C_{1}-2$.
As $c_{1}\ge0$ in the f\/irst place, we see that $m^{K,K_{1}}_{\mu}(\nu)=1$ if and only if (1)--(4) are satisf\/ied.
Suppose that $i>1$.
Without loss of generality we assume that $i=2$ and $j=3$ and
$\Xi=\{\epsilon_{1}\pm\epsilon_{2},\epsilon_{1}\pm\epsilon_{3}\}$.
Then
\begin{gather*}
m^{K,K_{1}}_{\mu}(\nu)
=p_{\Xi}(-c_{1}\epsilon_{1}+C_{2}\epsilon_{2}+C_{3}\epsilon_{3})-p_{\Xi}((c_{1}+2)\epsilon_{1}+C_{2}\epsilon_{2}+C_{3}\epsilon_{3}).
\end{gather*}
We have $-c_{1}\epsilon_{1}+C_{2}\epsilon_{2}+C_{3}\epsilon_{3}
=B_{2}(\epsilon_{2}-\epsilon_{1})+(C_{2}-B_{2})(\epsilon_{2}+\epsilon_{1})+B_{3}(\epsilon_{3}-\epsilon_{1})+(C_{3}-B_{3})(\epsilon_{3}+\epsilon_{1})$
from which it follows that
\begin{gather}
p_{\Xi}(-c_{1}\epsilon_{1}+C_{2}\epsilon_{2}+C_{3}\epsilon_{3})
\nonumber
\\
\qquad{}
=\#\left\{(B_{2},B_{3})\in\mathbb{N}^{2}:B_{2}\le C_{2},B_{3}\le C_{3}, B_{2}+B_{3}=\frac{1}{2}(C_{2}+C_{3}+c_{1})\right\}.
\label{eqn:line1}
\end{gather}
A~similar calculation shows
\begin{gather}
p_{\Xi}((c_{1}+2)\epsilon_{1}+C_{2}\epsilon_{2}+C_{3}\epsilon_{3})
\nonumber
\\
\qquad{}
=\#\left\{(B_{2},B_{3})\in\mathbb{N}^{2}:B_{2}\le C_{2},B_{3}\le C_{3}, B_{2}+B_{3}=\frac{1}{2}(C_{2}+C_{3}-c_{1}-2)\right\}.
\label{eqn:line2}
\end{gather}
The quantities~\eqref{eqn:line1}, \eqref{eqn:line2} can only be nonzero if $C_{2}+C_{3}-c_{1}$ is even.
In this case~\eqref{eqn:line1} is the number of integral points of the intersection of the line
$\ell=\{(B_{2},B_{3}):B_{2}+B_{3}=\frac{1}{2}(C_{2}+C_{3}+c_{1})\}$ with the rectangle $\{(B_{2},B_{3}):0\le B_{2}\le
C_{2}, 0\le B_{3}\le C_{3}\}$.
The quantity~\eqref{eqn:line2} is equal to the number of points of the intersection of the same rectangle with the line
$\ell+(1,0)$.
It follows that $m^{K,K_{1}}_{\mu}(\nu)=1$ if and only if $2\max(C_{2},C_{3})\le\frac{1}{2}(C_{2}+C_{3}+c_{1})\le
C_{2}+C_{3}$ which is equivalent to $|C_{2}-C_{3}|\le c_{1}\le C_{2}+C_{3}$.
\end{proof}

The branching from $K_{1}$ to~$M$ is equivalent to the branching $\mathrm{USp}(2)\times\mathrm{USp}(2)$ to the diagonal $\mathrm{USp}(2)$.
Given a~$K_{1}$ type $(c_{1}'\epsilon_{1}+c_{2}\epsilon_{2}+\dots+c_{n-1}\epsilon_{n-1}+c_{n}'\epsilon_{n})$,
the~$M$-types that occur upon restricting are of the form
\begin{gather*}
(c_{1}\epsilon_{1}+c_{2}\epsilon_{2}+\dots+c_{n-1}\epsilon_{n-1}+c_{1}\epsilon_{n}),
\qquad
c_{1}=\frac{|c_{1}'-c_{n}'|}{2},\frac{|c_{1}'-c_{n}'|}{2}+1,\ldots,\frac{c_{1}'+c_{n}'}{2}.
\end{gather*}
\begin{Corollary}
The matrix valued orthogonal polynomials of size $2\times2$ are obtained from $\mu=\omega_{1}$ or $\mu=\omega_{n-1}$.
\end{Corollary}
\begin{proof}
We have to do a~case by case investigation.
Let $\mu=x\omega_{i}+y\omega_{j}$ with $x$, $y$ non-negative integers,~$y$ positive and $i\le j$.
Consider the cases
\begin{itemize}
\itemsep=0pt
\item[$1)$] $i=j$,
\item[$2)$] $1\le i<n$, $j=n$,
\item[$3)$] $1\le i< j<n$.
\end{itemize}

Case~1.
If $i=n$ we have only one~$M$-type.
Indeed, this situation boiles down to restricting holomorphic or anti-holomorphic representation of
$\mathrm{SL}_{2}(\mathbb{C})$ to $\mathrm{SU}(2)$.
If $i=1$ then we have $\mu=(x+y)\omega_{1}$ and
\begin{gather*}
P_{K_{1}}^{+}(\mu)=\{s\epsilon_{1}+(x+y-s)\epsilon_{2}:s=0,\ldots,y+x\},
\end{gather*}
where $P_{K_{1}}^{+}(\mu)$ is the support of $m^{K,K_{1}}_{\mu}$.
If $i=n-1$ then $\mu=(x+y)\omega_{n-1}$ and
\begin{gather*}
P_{K_{1}}^{+}(\mu)=\{s\epsilon_{1}+(x+y)(\epsilon_{2}+\dots+\epsilon_{n-2})+s\epsilon_{n-1}:s=0,\ldots,y+x\}.
\end{gather*}
If $1<i<n-1$ then $\mu=(x+y)\omega_{i}$ then $(c_{i},c_{i+1})\in\{(0,0),(x+y,0),(x+y,x+y)\}$ lead to dif\/ferent
representations of $K_{1}$.

Case~2.
We have $\mu=x\omega_{i}+y\omega_{n}$ with $xy\ne0$.
If $1<i<n-1$ then we f\/ind three $K_{1}$ types that occur upon restriction.
If $i=1$ then we can choose $c_{2}=x$ and $c_{2}=0$.
The f\/irst choice implies $c_{1}=0$, the second $c_{1}=x$.
We f\/ind two $K_{1}$-types, one of them decomposes in at least two~$M$-types upon restriction.

Case 3.
If $x$, $y$ are both non zero, then we have at least three $K_{1}$-types.
Indeed, if either~$x$ or~$y$ is zero, then we are in Case 1.
Assume $xy\ne0$.
Then ${\mu,\mu-\epsilon_1-\epsilon_{i},\mu-\epsilon_1-\epsilon_{j}}$ are three $K_{1}$ types by Proposition~\ref{prop:
branching}.

We conclude that $\mu=\omega_{1}$ and $\mu=\omega_{n-1}$ are the only~$K$-types that give matrix valued orthogonal
polynomials of size $2\times2$.
\end{proof}

\subsection[The algorithm for a~$(2\times 2)$-case]{The algorithm for a~$\boldsymbol{(2\times 2)}$-case}\label{ss: 2x2-calc}

We consider $(G,K)=(\mathrm{USp}(2n),\mathrm{USp}(2n-2)\times\mathrm{USp}(2))$ and the irreducible
$\mathrm{USp}(2n-2)\times\mathrm{USp}(2)$-representations with highest weight~$\mu$ in a~two dimensional face.
The complexif\/ied Lie algebra of $\mathrm{USp}(2n)$ is denoted by $\mathfrak{sp}_{2n}(\mathbb{C})$.
We calculate the spherical function~$\phi$ and $\Psi^{\mu}_0$.
The realizations of the fundamental representations that are involved in this calculation are described in~\cite{Fulton-Harris}.
We employ the same notation as in~\cite{Fulton-Harris}.
In particular, the root vectors of $\epsilon_{i}-\epsilon_{j}$ and $2\epsilon_{k}$ are denoted by $X_{i,j}$ and $U_{k}$
respectively.

The spherical representation $\lambda_{\mathrm{sph}}=\varpi_{2}=\epsilon_{1}+\epsilon_{2}$ is realized in the kernel of
$\bigwedge^{2}V\to\mathbb{C}:v\wedge w\mapsto Q(v,w)$.
The kernel is of dimension $\binom{2n}{2}-1=n(2n-1)-1$.
The weight zero space $\bigwedge^{2}V(0)$ is spanned by $\{e_{k}\wedge e_{n+k}:k=1,\ldots,n\}$.
The weight zero space $V^{G}_{\varpi_{2}}(0)$ is spanned by the vectors $\{e_{k}\wedge e_{n+k}-e_{k+1}\wedge
e_{n+k+1}:k=1,\ldots,n-1\}$.
The vector in $V^{G}_{\varpi_{2}}$ that is f\/ixed under $\mathrm{USp}(2n-2)\times\mathrm{USp}(2)$ is of weight zero
because the groups~$G$ and~$K$ are both of semisimple rank~$n$.
We look for a~vector that is killed by $\mathfrak{sp}_{2n-2}(\mathbb{C})\oplus\mathfrak{sp}_{2}(\mathbb{C})$ and it is
easily seen that $e_{1}\wedge e_{n+1}+\ldots+e_{n-1}\wedge e_{2n-1}-(n-1)e_{n}\wedge e_{2n}$ is killed by this subalgebra.

The torus~$A$ has as its Lie algebra $\mathfrak{a}=\mathbb{R}\cdot (X_{1,n}-X_{n,1})$.
Put
\begin{gather*}
a_{t}=\left(
\begin{array}{@{}ccc|ccc@{}}
\cos t & 0 & \sin t& 0&0&0
\\
0 & I_{n-2} & 0 & 0&0&0
\\
-\sin t & 0 & \cos t& 0&0&0
\\
\hline
0 & 0 & 0 & \cos t &0&\sin t
\\
0 & 0 & 0 & 0 & I_{n-2}&0
\\
0 & 0 & 0 &-\sin t & 0 &\cos t
\end{array}
\right).
\end{gather*}
Then $A=\{a_{t}|t\in[0,2\pi)\}$ and $\langle v_{0},a_{t}v_{0}\rangle=n^2\cos^2(t)-n$ and it follows that
$\phi(a_{t})=(n\cos^{2}(t)-1)/(n-1)$.

Now we calculate the function $\Psi^{\omega_1}_0$ by means of the algorithm.
Upon identifying $V^{G}_{\varpi_{1}}\cong\mathbb{C}^{2n}$ we have $v_{\epsilon_{i}}=e_{i}$ and
$v_{-\epsilon_{i}}=e_{n+i}$.

\textbf{Step 1.}
We have $V_{\omega_{1}}|_{M}=\mathbb{C}^{2}\oplus\mathbb{C}^{2n-4}$, with $\mathbb{C}^{2}\cong\mathbb{C}
v_{\epsilon_{1}}\oplus\mathbb{C} v_{-\epsilon_{1}}$ and $\mathbb{C}^{2n-4}$ the standard representation of
$\mathrm{USp}(2n-4)$.
Hence $P_M^+(\omega_1)=\{\eta_1,\eta_2\}$.

\textbf{Step 2.}
The degree is zero, so we consider $B(\omega_1)=\{\varpi_{1},\varpi_{3}\}$

\textbf{Steps 4--10.}
We calculate $\Psi^{\omega_1}_{\varpi_1}(a_t)$.
The weight vectors of $V^{G}_{\varpi_{1}}$ are $\{v_{\pm\epsilon_{i}},1\le i\le n\}$.
The vector space $V^{K}_{\omega_{1}}$ is spanned by the weight vectors $v_{\pm\epsilon_{i}}$, $i=1,\ldots,n-1$ and the
embedding $V^{K}_{\omega_{1}}\to V^{G}_{\varpi_{1}}$ is clear.
We choose $v_{\varpi_1,\omega_1,\eta_1}=e_1$ and $v_{\varpi_1,\omega_1,\eta_2}=e_2$.
It follows that $\Psi^{\omega_{1}}_{\varpi_{1}}(a_{t})=(\cos(t),1)^T$.

\textbf{Steps 4--10:} we calculate $\Psi^{\omega_1}_{\varpi_3}(a_t)$.
Let $V=V^{G}_{\varpi_1}$ as above.
We realize $V^{G}_{\varpi_{3}}$ in the kernel of the map $\varphi_{3}:\bigwedge^{3}V\to V$ given~by
\begin{gather*}
v_{k_{1}}\wedge v_{k_{2}}\wedge
v_{k_{3}}\mapsto\sum\limits_{i<j}Q(v_{i},v_{j})(-1)^{i+j-1}\cdots\wedge\widehat{v_{i}}\wedge\cdots\wedge\widehat{v_{j}}\wedge\cdots.
\end{gather*}
In $\bigwedge^{3}V$, the weight vectors of weight $\pm\epsilon_{i}\pm\epsilon_{j}\pm\epsilon_{k}$ are
$v_{\pm\epsilon_{i}}\wedge v_{\pm\epsilon_{j}}\wedge v_{\pm\epsilon_{k}}$.
There are $8\binom{n}{3}$ of them.
The weight vectors of weight $\pm\epsilon_{i}$ are $v_{\epsilon_{j}}\wedge v_{-\epsilon_{j}}\wedge v_{\pm\epsilon_{i}}$.
There are $2n(n-1)$ of them.
The kernel of this map is of dimension $\binom{2n}{3}-2n=8\binom{n}{3}+2(n-2)n$, the multiplicities of the short roots
are $n-2$ and there are $2n$ of them.
It follows that the restricted map $\bigwedge^{3}V(\pm\epsilon_{i})\to V$ has an $(n-2)$-dimensional kernel.
Since
\begin{gather*}
\varphi^{3}(v_{\epsilon_{j}}\wedge v_{-\epsilon_{j}}\wedge
v_{\pm\epsilon_{i}})=(-1)^{n-1}Q(v_{\epsilon_{j}},v_{-\epsilon_{j}})v_{\pm\epsilon_{i}},
\end{gather*}
we have
\begin{gather*}
\ker\big(\varphi^{3}\big|_{\bigwedge^{3}V(\epsilon_{i})}\big)
=\left\{\sum\limits_{j\ne i}a_{j}v_{\epsilon_{j}}\wedge v_{-\epsilon_{j}}\wedge v_{\epsilon_{i}}\Bigg|\sum\limits_{j\ne i}a_{j}=0\right\}.
\end{gather*}
In order to calculate the embedding $V^{K}_{\omega_{1}}$ into $V^{G}_{\varpi_{3}}$ we need to calculate
\begin{gather*}
\bigcap_{\alpha\in R_{K}^{+}}\ker\Big(e_{\alpha}\big|_{\ker(\varphi^{3}|_{\bigwedge^{3}V(\epsilon_{1})})}\Big).
\end{gather*}
It is suf\/f\/icient to calculate the kernels of the root vectors for the simple roots.
Note that $X_{1,2}$ already acts as zero.
For the others we have
\begin{gather*}
X_{k,k+1}\left(\sum\limits_{j=2}^{n}a_{j}v_{\epsilon_{j}}\wedge v_{-\epsilon_{j}}\wedge v_{\epsilon_{1}}\right)
=(a_{k}-a_{k+1})v_{\epsilon_{k}}\wedge v_{-\epsilon_{k+1}}\wedge v_{\epsilon_{1}}
\end{gather*}
for $k=2,\ldots,n-2$ and
\begin{gather*}
U_{k}\left(\sum\limits_{j=2}^{n}a_{j}v_{\epsilon_{j}}\wedge v_{-\epsilon_{j}}\wedge v_{\epsilon_{1}}\right)=0
\end{gather*}
for $k=n-1,n$.
Hence $a_{k}=a_{k+1}$ for $k=2,\ldots,n-2$ and
\begin{gather*}
\bigcap_{\alpha\in R_{K}^{+}}\ker\Big(e_{\alpha}\big|_{\ker(\varphi^{3}|_{\bigwedge^{3}V(\epsilon_{1})})}\Big)=
\mathbb{C}\left(\sum\limits_{j=2}^{n-1}v_{\epsilon_{j}}\wedge v_{-\epsilon_{j}}\wedge
v_{\epsilon_{1}}-(n-2)(v_{\epsilon_{n}}\wedge v_{-\epsilon_{n}}\wedge v_{\epsilon_{1}})\right).
\end{gather*}
It follows that $V^{K}_{\omega_{1}}\to V^{G}_{\varpi_{3}}$ is determined~by
\begin{gather*}
v_{\epsilon_{1}}\mapsto w:=\sum\limits_{j=2}^{n-1}v_{\epsilon_{j}}\wedge v_{-\epsilon_{j}}\wedge
v_{\epsilon_{1}}-(n-2)(v_{\epsilon_{n}}\wedge v_{-\epsilon_{n}}\wedge v_{\epsilon_{1}}).
\end{gather*}
The highest weight vector of the~$M$-type $\eta_2$ in $V^{K}_{\omega_1}$ is $X_{2,1}v_{\epsilon_{1}}$ which maps to
\begin{gather*}
X_{2,1}w=X_{2,1}\left(\sum\limits_{j=2}^{n-1}v_{\epsilon_{j}}\wedge v_{-\epsilon_{j}}\wedge
v_{\epsilon_{1}}-(n-2)(v_{\epsilon_{n}}\wedge v_{-\epsilon_{n}}\wedge v_{\epsilon_{1}})\right)
\\
\phantom{X_{2,1}w}
=v_{\epsilon_{1}}\wedge v_{-\epsilon_{1}}\wedge v_{\epsilon_{2}}+\sum\limits_{j=3}^{n-1}v_{\epsilon_{j}}\wedge
v_{-\epsilon_{j}}\wedge v_{\epsilon_{2}}-(n-2)(v_{\epsilon_{n}}\wedge v_{-\epsilon_{n}}\wedge v_{\epsilon_{2}})\in
V_{\varpi_3}.
\end{gather*}
Hence $(a_{t}w,w)=(n-2)\cos(t)+(n-2)^{2}\cos(t)$ and $(a_{t}X_{2,1}w,X_{2,1}w)=(n-1)^{2}\cos^{2}(t)-(n-1)$, so that
$\Psi^{\omega_1}_{\varpi_3}(a_t)=(\cos(t),((n-1)\cos^2(t)-1)/(n-2))^T$.
We f\/ind
\begin{gather*}
\Psi^{\omega_{1}}_{0}(a_{t})=
\begin{pmatrix}
\cos(t) & \cos(t)
\\
1 & ((n-1)\cos^{2}(t)-1)/(n-2)
\end{pmatrix}
.
\end{gather*}
We resume this discussion in the following statements.
\begin{Proposition}
\label{prop:Psi_0_w1}
In the variable $x=\cos^2(t)$, we have
\begin{gather*}
\widetilde \Psi_0^{\omega_1}(x)=
\begin{pmatrix}
\sqrt{x} & \sqrt{x}
\vspace{1mm}\\
1 & \dfrac{x(n-1)-1}{n-2}
\end{pmatrix}.
\end{gather*}
\end{Proposition}

\begin{Theorem}
\label{thm:example_2x2_SPn_1}
The multiplicity free triple $(\mathrm{USp}(2n),\mathrm{USp}(2n-2)\times\mathrm{USp}(2),\omega_1)$ gives rise to the
MVCPs $(W^{\omega_1}_{n},D^{\omega_1}_{n})$ given~by
\begin{gather*}
W^{\omega_1}_{n}(x)=x(1-x)^{2n-3}
\begin{pmatrix}
2x+2n-4 & 2xn-2
\vspace{1mm}\\
2xn-2 & 2\dfrac {((n-1)^2x^2-nx+1)}{n-2}
\end{pmatrix},
\\
\widetilde D^{\omega_1}=x(x-1)\partial_x^2+\big(-2-2\widetilde{S}+2x(n-\widetilde R)\big)\partial_x+\Lambda_0/2,
\end{gather*}
with
\begin{gather*}
\widetilde{R}=
\begin{pmatrix}
-\dfrac12 & \dfrac12
\vspace{1mm}\\
0 & -1
\end{pmatrix},
\qquad
\widetilde S = \frac{1}{2(n-1)}
\begin{pmatrix}
1 & 1
\\
(n-2) & (n-2)
\end{pmatrix},
\qquad
\Lambda_0=
\begin{pmatrix}
0 & 0
\\
0 & 4(n-1)
\end{pmatrix}.
\end{gather*}
\end{Theorem}
\begin{proof}
We follow the steps in Section~\ref{subsection:calculation data}.
From the explicit expression of $\widetilde \Psi_0^{\omega_1}(x)$ in Proposition~\ref{prop:Psi_0_w1}, we compute
$x(1-x)(\widetilde \Psi_0^\mu(x))^{-1}\partial_{x}\widetilde \Psi_0^ \mu(x)$ which gives the matrices $\widetilde R$, $\widetilde S$.
We determine the expression for the eigenvalue $\Lambda_0$ by calculating how the Casimir operator acts on the
representation spaces at hand.
Remark~\ref{remark:list_of_operators} and Table~\ref{table2} give the expression of the dif\/ferential ope\-ra\-tor~$\widetilde D^{\omega_1}$.
Observe that $\Lambda_0$ is normalized so that the $(1,1)$-entry is zero by adding a~multiple of the identity matrix.
\end{proof}

\begin{Remark}
\label{rmk:TiraoMVHF}
Let us consider a~dif\/ferential operator of the form
\begin{gather}
\label{eq:Tirao-hypergeometric_equation}
z(1-z)F''(z) + (C-zU)F'(z)  -  VF(z)=0,
\qquad
z\in\mathbb{C},
\end{gather}
where~$C$,~$U$ and~$V$ are $N\times N$ matrices and $F\colon \mathbb{C}\to \mathbb{C}^N$ is a~(column-)vector-valued
function which is twice dif\/ferentiable.
It is shown by Tirao~\cite{Tirao-PNAS} that if the eigenvalues of~$C$ are not in~$-\mathbb{N}$, then the matrix-valued
hypergeometric function~${}_2H_1$ def\/ined as the power series
\begin{gather*}
\tHe{U, V}{C}{z}   =  \sum\limits_{i=0}^\infty \frac{z^i}{i!}   [C,U,V]_i,
\\
[C,U,V]_0=1,
\qquad
[C,U,V]_{i+1}=(C+i)^{-1}\bigl(i^2+i(U-1)+V\bigr) [C,U,V]_i
\end{gather*}
converges for $|z|<1$ in $M_d(\mathbb{C})$.
Moreover, for $F_0\in \mathbb{C}^N$ the (column-)vector-valued function
\begin{gather}
\label{eq:F(z)_hyperg}
F(z) \, = \, \tHe{U, V}{C}{z}\, F_0
\end{gather}
is a~solution to~\eqref{eq:Tirao-hypergeometric_equation} which is analytic for $|z|<1$, and any analytic (on $|z|<1$)
solution to~\eqref{eq:Tirao-hypergeometric_equation} is of this form.
\end{Remark}

In the following Corollary we write the monic MVOPs with respect to $W^{\omega_1}$ in terms of the matrix-valued
hypergeometric functions ${}_2H_1$.
\begin{Corollary}
\label{cor:sf_as_Hypergeometric}
The unique sequence of monic MVOP $\{P_d\}_{d\geq0}$ with respect to $W^{\varpi_1}$ is given~by
\begin{gather*}
\big(P_d(x)M^{-1}\big)_{i,j} =\!\left(\!\tHe{2(n-\widetilde R), \Lambda_0/2-\lambda_d(i)}{2\widetilde S+2}{x}\! F_d(i)\!\right)_j\!, \!\qquad M=
\begin{pmatrix}
1 & \dfrac{d}{d+2n-2}
\vspace{1mm}\\
0 & 1
\end{pmatrix}\!,
\end{gather*}
where $\lambda_d(1)=d(d+2n)$, $\lambda_d(2)=d(d+2n+1)+2n-2$ and
\begin{gather*}
F_d(i)=d!
\big[2\widetilde S+2, 2(n-\widetilde R), \Lambda_0/2-\lambda_d(i)\big]_d^{-1} e_i.
\end{gather*}
Here $e_i$ is the standard basis vector.
\end{Corollary}
\begin{proof}
Since the eigenvalues of $2+2\widetilde S$ are not in $-\mathbb{N}$, we can apply Remark~\ref{rmk:TiraoMVHF}.
By looking at the leading coef\/f\/icient, it is easy to see that the eigenvalue of $\widetilde D^{\omega_1}$ for $P_d$ is
given by the upper triangular matrix $\widetilde \Lambda_d=d(d-1)+d2(n-\widetilde R)+\Lambda_0/2$.
It is readily seen that the polynomial $P_d(x)M^{-1}$, $d\in \mathbb{N}_0$, is an eigenfunction of $\widetilde
D^{\omega_1}$ with eigenvalue
\begin{gather*}
\widetilde \Lambda_d^M=M\widetilde \Lambda_d M^{-1} =
\begin{pmatrix}
\lambda_d(1) & 0
\\
0 & \lambda_d(2)
\end{pmatrix}
=
\begin{pmatrix}
d(d+2n) & 0
\\
0 & d(d+2n+1)+2n-2
\end{pmatrix}.
\end{gather*}
Therefore the~$i$-th column of $P_dM^{-1}$, $i=1,2$, is a~polynomial solution of the equation
\begin{gather}
\label{eq:hyperg_example1_spn}
x(1-x)\partial_x^2+\big(2+2\widetilde{S}-2x(n-\widetilde R)\big)\partial_x-(\Lambda_0/2-\lambda_d(i))=0.
\end{gather}
It can be easily verif\/ied that $\lambda_d(i)=\lambda_{d'}(i')$ if and only if $d=d'$ and $i=i'$, see for
instance~\cite[Lemma~2.2]{PachR}.
This implies that there is a~unique (up to a~scalar multiple) polynomial solution of~\eqref{eq:hyperg_example1_spn}, see
also~\cite[Theorem~4.5]{KvPR2} and the discussion above.
The explicit expression for the columns of $P_d(x)M^{-1}$ follows from~\eqref{eq:F(z)_hyperg}.
\end{proof}

\subsection[The second family of $2\times2$ MVOP]{The second family of $\boldsymbol{2\times2}$ MVOP}\label{section:second_2x2_example}

In this subsection we use our implementation in \verb|GAP|~\cite{vPR} to obtain a~one-parameter sequence of MVOPs.
We conjecture that it is the one associated to the spherical functions of type $\mu_n=\varpi_{n-1}$ for the pair
$(G,K)=(\mathrm{USp}(2n),\mathrm{USp}(2n-2)\times\mathrm{USp}(2))$.
As we were not able to calculate $\Psi_{0}^{\mu}$ by hand, we computed $\widetilde \Psi^{\mu_n}_0$ for small values
of~$n$ in \verb|GAP| and made an ansatz for its general expression.

\begin{Conjecture}
For all $n\geq 3$, we have
\begin{gather*}
\widetilde \Psi^{\mu_n}_0(x) =
\begin{pmatrix}
\dfrac{x+1}{2} & \dfrac{(n+1)x-2}{n-1}
\vspace{1mm}\\
\sqrt{x} & \dfrac{\sqrt{x}((n+3)x+n-5)}{2(n-1)}
\end{pmatrix}.
\end{gather*}
\end{Conjecture}

The matrix $\widetilde \Psi^{\mu_n}_0$ is the building block of the weight matrix $W^{\mu_n}$.
If this function is indeed the desired function $\Psi_{0}^{\mu}$ that comes from the representation theory, then the
construction to make a~MVCP out of it should work.
To collect the matrices that appear as coef\/f\/icients in the dif\/ferential operator, we compute explicitly the f\/irst-order
dif\/ferential equation.
\begin{Lemma}
\label{lem:RS_2x2}
The matrix $\widetilde \Psi^{\mu_n}_0$ satisfies~\eqref{eq:first_order_equation_x} with
\begin{gather*}
\widetilde S=
\begin{pmatrix}
\dfrac{4}{n+3} & \dfrac{2n-10}{(n-1)(n+3)}
\vspace{1mm}\\
\dfrac{n-1}{n+3} & \dfrac{n-5}{2(n+3)}
\end{pmatrix},
\qquad
\widetilde R=
\begin{pmatrix}
-1 & \dfrac{(n+1)}{(n-1)}
\vspace{1mm}\\
0 & -\dfrac{3}{2}
\end{pmatrix}.
\end{gather*}
\end{Lemma}
\begin{proof}
The lemma follows by computing explicitly
\begin{gather*}
x(1-x)\big(\widetilde \Psi^{\mu_n}_0(x)\big)^{-1} \frac{d}{dx}\widetilde \Psi^{\mu_n}_0(x).
\end{gather*}
Then $\widetilde S$ is minus the coef\/f\/icient of degree $0$ in~$x$ and $\widetilde R$ is the coef\/f\/icient of degree 1
in~$x$.
\end{proof}

\begin{Theorem}
\label{thm:second_example_2x2}
Let $W^{\mu_n}$ and $\widetilde D^{\mu_n}$ be defined~by
\begin{gather*}
W^{\mu_n}(x)=(1-x)^{2n-3}x\widetilde \Psi^{\mu}_{0}(x)^{*}T^{\mu_{n}}\widetilde \Psi^{\mu}_{0}(x),
\qquad
T^{\mu_{n}}=
\begin{pmatrix}
\dfrac{2}{n+1} & 0
\vspace{1mm}\\
0 & \dfrac{(2n-2)}{(n-1)(n-2)}
\end{pmatrix}.
\\
\widetilde D^{\mu_{n}}=x(1-x)\partial_x^2 + \big[2+2\widetilde S -2x(n-\widetilde R)\big]\partial_x -\Lambda_0/2,
\qquad
\Lambda_0=
\begin{pmatrix}
0 & 0
\\
0 & 2n+6
\end{pmatrix},
\end{gather*}
with $\widetilde R$ and $\widetilde S$ given in Lemma~{\rm \ref{lem:RS_2x2}}.
Then $(W^{\mu_n},\widetilde D^{\mu_n})$ is a~MVCP.
\end{Theorem}
\begin{proof}
The proof is analogous to that of Theorem~\eqref{thm:example_2x2_SPn_1}.
The matrix $T^{\mu}$, or rather a~multiple of it, is computed by means of Weyl's dimension formula.
Finally we use the explicit expressions of $(W^{\mu_n},\widetilde{D}^{\mu_n})$ to verify the symmetry
equations~\eqref{eq:symmetry_equations}.
\end{proof}

\begin{Remark}
The MVCP that we obtained in Theorem~\ref{thm:second_example_2x2} is only by conjecture the MVCP associated to the
indicated multiplicity free triple.
The statement is somewhat misleading because we persisted in using the same notation for the weight and the dif\/ferential operator.
However, the fact that the construction works for the conjectured expression of $\Psi_{0}^{\mu}$, convinces us that it is the right one.
Indeed, computer experiments show that it is unlikely that the wrong expression of $\Psi_{0}^{\mu}$, yield a~MVCP after all.
\end{Remark}

\begin{Remark}
Observe that, as in the previous example, the eigenvalues of $\widetilde S$ are $0$ and $1/2$.
This implies that the eigenvalues of $2+2\widetilde S$ are $2$ and $3$ so that the polynomial eigenfunctions of the
dif\/ferential operator in Theorem~\ref{thm:second_example_2x2} can be written in terms of matrix valued hypergeometric
functions~\cite{Tirao-PNAS}.
\end{Remark}

\subsection[A family of $3\times3$ MVOP]{A family of $\boldsymbol{3\times3}$ MVOP}

In this subsection we f\/ix the Gelfand pair $(G,K)=(\mathrm{USp}(6), \mathrm{USp}(4)\times \mathrm{USp}(2))$ and we study
the one-parameter family of~$K$-types $\mu_j=\varpi_1+j\varpi_3$.
For small values of $j\in \mathbb{N}$ we can use our implementation in GAP to compute $\widetilde \Psi^{\mu_j}_0$.
We use this information to make the following ansatz of $\widetilde \Psi^{\mu_j}_0$ for all~$j$: For all
$j\in\mathbb{N}$, we have
\begin{gather*}
\widetilde \Psi^{\mu_j}_0 =
\begin{pmatrix}
\dfrac{x^{\frac{(j-1)}{2}}((j+1)x-1)}{j} & x^{\frac{(j-1)}{2}} & \dfrac{x^{\frac{(j-1)}{2}}(1+4x)}{5}
\vspace{1mm}\\
x^{\frac{(j-1)}{2}} & x^{\frac{(j-1)}{2}} (-j+x(j+1)) & x^{\frac{(j-1)}{2}}\dfrac{(-1+x(j+5))}{5}
\vspace{1mm}\\
x^{\frac{j}{2}} & x^{\frac{j}{2}} & x^{\frac{j}{2}}\dfrac{(-(2j+5)+ 2x(j+5))}{5}
\end{pmatrix}.
\end{gather*}
As in the previous example, we obtain the matrices $\widetilde R$ and $\widetilde S$ in the dif\/ferential
equation~\eqref{eq:Dmu} for $\widetilde \Psi^{\mu_j}_0$
\begin{gather*}
\widetilde S =
\begin{pmatrix}
-\dfrac{(j^3+7j^2+7j-5)}{2(j+1)(j+5)} & -\dfrac{j(j+4)}{(j+1)(j+5)} & -\dfrac{j(3j+10)}{5(j+1)(j+5)}
\vspace{1mm}\\
-\dfrac{3}{2(j+2)} & -\dfrac{(2j^2+1)}{4(j+1)} & \dfrac{(2j-5)}{20(j+1)}
\vspace{1mm}\\
-\dfrac{5}{2(2j+10)} & -\dfrac{5}{2(2j+10)} & -\dfrac{(2j^2+10j+5)}{4(j+5)}
\end{pmatrix},
\\
\widetilde R =
\begin{pmatrix}
-\dfrac{(j+1)}{2} & 0 & \dfrac{2j}{5(j+1)}
\vspace{1mm}\\
0 & -\dfrac{(j+1)}{2} & \dfrac{(j+5)}{10(j+1)}
\vspace{1mm}\\
0 & 0 & -\dfrac{(j+2)}{2}
\end{pmatrix}
\end{gather*}
Now we proceed as in Section~\ref{subsection:calculation data} to construct the MVCP associated to~$\mu_j$.
\begin{Theorem}
Let $W^{\mu_n}$ and $D^{\mu_n}$ be defined~by
\begin{gather}
W^{\mu_n}=(1-x)^3x\widetilde \Psi^{\mu_j}_{0}(x)^{*}T^{\mu_j}\widetilde \Psi^{\mu_j}_{0}(x),
\nonumber
\\
\label{eq:diff_operator_3x3}
\widetilde D^{\mu}=x(1-x)\partial_{x}^2+\big[2+2\widetilde S - 2x(n-\widetilde R)\big]
\partial_{x}-\frac{\Lambda_0}{2},
\end{gather}
where
\begin{gather*}
T^{\mu_j}=
\begin{pmatrix}
i & 0 & 0
\\
0 & i+2 & 0
\\
0 & 0 & 2i+2
\end{pmatrix},
\qquad
\Lambda_0=
\begin{pmatrix}
0 & 0 & 0
\\
0 & i+1 & 0
\\
0 & 0 & i+5
\end{pmatrix}.
\end{gather*}
Then $(W^{\mu_n}, \widetilde D^{\mu_n})$ is a~MVCP.
\end{Theorem}
\begin{Remark}
The eigenvalues of $\widetilde S$ are $(i+1)/2$, $i/2$, $(i-1)/2$ so that $2+\widetilde S$ has always positive eigenvalues.
This ensures that the hypergeometric operator~\eqref{eq:diff_operator_3x3} f\/its into the theory developed in~\cite{Tirao-PNAS}.
However, further analysis is required to f\/ind an expression of the orthogonal polynomial with respect to $W^{\mu_n}$ as
matrix hypergeometric functions.
\end{Remark}

\section[Classif\/ication of all $(2\times2)$-cases]{Classif\/ication of all $\boldsymbol{(2\times2)}$-cases}\label{2x2class}

In this section we give the other MVCPs of size $2\times2$ from Table~\ref{table mfs compact}.
First we classify all the possible $2\times2$ cases, then we indicate the general approach of how we calculate the
matrix coef\/f\/icients and then we give explicit expressions for the MVCPs.
The fundamental weights of~$G$ are denoted by $\varpi_{i}$, those of~$K$ by $\omega_{i}$ and those of~$M$ by $\eta_{i}$.
We employ the standard def\/initions and conventions for weights according to~\cite{Knapp}.
The embeddings $K\to G$ are those as in~\cite{HvP}.

\begin{Theorem}
Let $(G,K,F)$ be a~multiplicity free system of Table~\ref{table mfs compact} and let $\mu\in F$.
Then $(W^{\mu},D^{\mu})$ is of size $2\times 2$ if and only if $(G,K,\mu)$ is one of the following:
\begin{itemize}
\itemsep=0pt
\item[$a)$] $(G,K)=(\mathrm{SU}(n+1),\mathrm{U}(n))$ and $\mu=\varpi_{i}+m\varpi_{n}$ with $i\in\{1,\ldots,n-1\}$
and $m\in\mathbb{Z}$;
\item[$b)$] $(G,K)=(\mathrm{SO}(2n+1),\mathrm{SO}(2n))$ and $\mu=\omega_{i}$ with $i\in\{1,\ldots,n-2\}$;
\item[$d)$] $(G,K)=(\mathrm{SO}(2n),\mathrm{SO}(2n-1))$ and $\mu=\omega_{i}$ with $i\in\{1,\ldots,n-1\}$;
\item[$c)$] $(G,K)=(\mathrm{USp}(2n),\mathrm{USp}(2n-2)\times\mathrm{USp}(2))$ and
$\mu\in\{\omega_{1},\omega_{n-1}\}$;
\item[$g)$] $(G,K)=(\mathrm{G}_{2},\mathrm{SU}(3))$ and $\mu=\omega_{1}$ or $\mu=\omega_{2}$;
\item[$f)$] $(G,K)=(\mathrm{F}_{4},\mathrm{Spin}(9))$ and $\mu=\omega_{1}$.
\end{itemize}
\end{Theorem}
\begin{proof}
We have to prove that the indicated irreducible~$K$-representations decompose into two irreducible~$M$-representations
upon restriction to~$M$.
In all cases apart from Case~f this follows from the branching rules $K\downarrow M$ described in~\cite{Camporesi2}
(Cases~a,~b and~d), Section~\ref{section:SPcase} (Case~c) and~\cite{HvP} (Case~g).
To exclude all the other cases for the pairs $(G,K)$ in Cases~a,~b and~d  we refer to the same branching rules.
These are all given by interlacing conditions and the statement is easily checked.
Roughly speaking, the~$K$-weight is given by a~string of ordered integers and as soon as we make more than one jump, or
a~jump larger than one, we will f\/ind more than two~$M$-weights that interlace.
Case~c is dealt with in Section~\ref{section:SPcase}.
For Case~g we have to check the branching $\mathrm{SU}(3)\downarrow\mathrm{SU}(2)$ which can also be done in stages,
$\mathrm{SU}(3)\downarrow\mathrm{U}(2)\downarrow\mathrm{SU}(2)$.
On the f\/irst step, branching Case~a applies.

For Case~f we use the branching rules $\mathrm{Spin}(9)\downarrow\mathrm{Spin}(7)$ as described
in~\cite{Baldoni-Silva,HvP} to see that
$|P^{+}_{M}(a\omega_{1}+b\omega_{2})|=\sum\limits_{c=0}^{a}\sum\limits_{d=0}^{b}\sum\limits_{e=0}^{d}1$,
$|P_{M}^{+}(a\omega_{3})|=\sum\limits_{b=0}^{a}\sum\limits_{c=0}^{b}\sum\limits_{d=0}^{a-b}1$ and
$|P^{+}_{M}(a\omega_{4})|=2a+2$ if $a>0$ (and one if $a=0$).
It follows that $\mu=\omega_{1}$ yields the only $(2\times 2)$-case for the pair $(\mathrm{F}_{4},\mathrm{Spin}(9))$.

To exclude the pair $(\mathrm{Spin}(7),\mathrm{G}_{2})$, we use the branching rule for
$\mathrm{G}_{2}\downarrow\mathrm{SU}(3)$.
The $\mathrm{SU}(3)$-types that occur are parametrized by the integral points of a~triangle with integral vertices,
i.e.~there are no $2\times2$ cases.
\end{proof}

In order to calculate the MVCPs we have to determine the function $\widetilde \Psi^{\mu}_{0}$ and the fundamental
spherical function~$\phi$.
The latter is already used in the proof of Lemma~\ref{lemma: mM}.
To calculate $\widetilde \Psi_{0}^{\mu}$ we introduce the following notation.
Let $B(\mu)=\{\lambda,\lambda'\}$ and let $P^{+}_{M}(\mu)=\{\nu,\nu'\}$, where $\lambda=\lambda(0,\nu)$ and
$\lambda'=\lambda(0,\nu')$ in the notation of Theorem~\ref{thm:parametrization mu-well}.
Let $V^{M}_{\nu}$, $V^{K}_{\mu}$, $V^{G}_{\lambda}$ denote the representations spaces and consider the following commutative diagram
\begin{gather}
\begin{split}
&\begin{xy}
\xymatrix{
   & V^{G}_{\lambda}\\
V^{M}_{\nu}\oplus V^{M}_{\nu'}\ar[r]^{i} \ar[ur]^f\ar[dr]_{f'}& V^{K}_{\mu}\ar[u]_{\gamma}\ar[d]^{\gamma'}\\
 & V^{G}_{\lambda'}
}
\end{xy}
\end{split}
\label{diagram}
\end{gather}

The maps $f$, $f'$ are~$M$-equivariant,~$i$ is an~$M$-isomorphism and $\gamma$, $\gamma'$ are the embeddings of
the~$K$-representation spaces and are thus~$K$-equivariant.
Let $w_{\nu}$, $w_{\nu'}$ denote the highest weight vectors of the~$M$-representations.
We need to calculate the four vectors $x_{1,1}=f(w_{\nu})$, $x_{2,1}=f(w_{\nu'})$, $x_{1,2}=f'(w_{\nu})$ and
$x_{2,2}=f'(w_{\nu'})$ and subsequently the inner products $\langle x_{i,j},a_{t}x_{i,j}\rangle$.
Together with information from Table~\ref{table2} we calculate the data $(\widetilde
\Psi^{\mu}_{0},\phi,T^{\mu},\lambda_{1},\Lambda_{0})$.
Then we follow the steps in Section~\ref{subsection:calculation data} to obtain explicit expressions for the
corresponding pair $(W^{\mu}, \widetilde D^{\mu})$.

In the remainder of this section we provide expressions for almost all the functions $\widetilde \Psi_0^\mu$ of size
$2\times 2$ that can be obtained from Table~\ref{table mfs compact} using our method (note that Case~c is
dealt with in Section~\ref{section:SPcase}).
See Remark~\ref{remark:real parameters} for further comments on the parameters.

Case~a. $(G,K)=(\mathrm{SU}(n+1),\mathrm{U}(n))$.
In this case we need to distinguish between two cases.
The f\/irst case has been already investigated in a~series of papers ending with~\cite{PT}.

Case~a1. We have a~two-parameter family of of classical pairs corresponding to the~$K$-types
$\mu=\varpi_{i}+m\varpi_{n}$, where $i\in\{1,\ldots,n-1\}$ and $m\in\mathbb{Z}_{\geq 0}$.
In this case the bottom of the~$\mu$-well is given by $B(\mu)=\{\mu,\mu+\epsilon_{i+1}-\epsilon_{n+1}\}$.
The function $\widetilde \Psi_0^\mu$ can be obtained by applying our algorithm or from~\cite[Section~5.2]{PT}
\begin{gather*}
\widetilde \Psi_0^{\mu}(x)=x^{\frac{m}{2}}
\begin{pmatrix}
\sqrt{x} & \sqrt{x}
\vspace{1mm}\\
1 & \dfrac{(m+1)-x(m+n-i+1)}{i-n}
\end{pmatrix}.
\end{gather*}

Case~a2. This case was not considered in the literature before as far as we know.
We take the~$K$-types $\mu=\varpi_{i}+m\varpi_{n}$ but now we let $m\in\mathbb{Z}_{< 0}$.
Note that although all ingredients needed for this case are already given in~\cite{PT}, the weight matrix constructed
there does not have f\/inite moments for negative values of~$m$, thus it is excluded~\cite[Section 6]{PT}.
Our case is essentially a~conjugation of the case considered in~\cite{PT} and does not have any problem of
integrability.
The function $\Psi_0^\mu$ is given~by
\begin{gather*}
\widetilde \Psi_0^{\mu}(x)=x^{-\frac{(m+1)}{2}}
\begin{pmatrix}
\dfrac{(m+(i-m)x))}{i} & 1
\vspace{1mm}\\
x^{\frac12} & x^{\frac12}
\end{pmatrix}.
\end{gather*}
Cases~b and~d.
We only treat $(G,K)=(\mathrm{SO}(2n+1),\mathrm{SO}(2n))$ and $\mu=\omega_{i}$ with $i\in\{1,\ldots,n-2\}$.
This case is essentially the same as the one in Section~\ref{section:SPcase} but simpler.
We use the notation of~\cite[Chapter~19]{Fulton-Harris}.
The~$K$ representation is of highest weight $\omega_{i}$ for indicated~$i$ is realized in $\bigwedge^{i}\mathbb{C}^{2n}$
and this space decomposes as $\bigwedge^{i-1}\mathbb{C}^{2n-1}\oplus\bigwedge^{i}\mathbb{C}^{2n-1}$
as~$M$-representation, where $M\cong\mathrm{SO}(2n-1)$.
It follows that $P^{+}_{M}(\omega_{i})=\{\eta_{i-1},\eta_{i}\}$.
On the other hand, $B(\omega_{i})=\{\varpi_{i},\varpi_{i+1}\}$ and $\lambda(0,\eta_{i})=\varpi_{i+1}$.
Everything is now explicit and the maps $f$, $f'$ of the commutative diagram~\eqref{diagram} are easily calculated on
highest weight vectors.
We have $A=\{a_{t}=\exp(tH_{A}):t\in\mathbb{R}\}$, where $H_{A}=X_{1,n}-X_{n,1}$.  We f\/ind
\begin{gather*}
\Psi^{\omega_{i}}_{0}(a_{t})= \left(
\begin{matrix} \cos(t)&1
\\
1& \cos(t)
\end{matrix}
\right).
\end{gather*}

Case~g1.
$(G,K)=(\mathrm{G}_{2},\mathrm{SU}(3))$ and $\mu=\omega_{1}$.
Let $\mu=\omega_{1}$.
We have $M\cong\mathrm{SU}(2)$.
Following~\cite{HvP} we f\/ind $P^{+}_{M}(\mu)=\{0,\eta_1\}$ and $B(\mu)=\{\varpi_{1},\varpi_{2}\}$, with
$\lambda(0,0)=\varpi_{1}$ and $\lambda(0,\eta_1)=\varpi_2$.
Fol\-lowing~\cite[Paragraph~2.2.7]{vP} we see that the Lie algebra of~$A$ is generated by $E_{\varpi_{1}}-E_{-\varpi_{1}}$.
It follows, as in all cases, that the calculations of the restricted matrix coef\/f\/icients in~$X_{\varpi_{1}}$ and~$X_{\varpi_{2}}$ are really $\mathrm{SU}(2)$-calculations.
From the weight diagrams (see, e.g.,~\cite[Chapter~22]{Fulton-Harris}) we read:
\begin{gather*}
\Psi^{\omega_{i}}_{0}(a_{t}) = \left(
\begin{matrix} m^{1}_{1,1}(b_{2t}) & m^{1}_{1,1}(b_{2t})
\vspace{1mm}\\
m^{1/2}_{1/2,1/2}(b_{2t}) & m^{3/2}_{1/2,1/2}(b_{2t})
\end{matrix}
\right)
 =\left(
\begin{matrix} \cos^{2}(t)&\cos^{2}(t)
\vspace{1mm}\\
\cos(t)&\cos^3(t)-2\cos(t)\sin^{2}(t)
\end{matrix}
\right),
\end{gather*}
where the $m^{\ell}_{m,n}$ denote the matrix coef\/f\/icients of $\mathrm{SU}(2)$ and where $b_{2t}$ is the rotation over an
angle~$t$, see~\cite{Koornwinder}.

The Cases~g2 and~f are left: Case~g2 is similar to Case~g1, and Case~f is too complicated to do~by
hand.
The involved $F_{4}$ representations are of dimension~52 and~1274.
We will handle it as soon as our \verb|GAP| implementation is improved to deal with general cases.

\subsection*{Acknowledgements}

The research for this paper was partly conducted when the f\/irst author visited the University of C{\'o}rdoba in August
and September 2012.
We would like to thank the Mathematics Departments of the Universities of Nijmegen and C{\'o}rdoba for their generous
supports that made this visit possible.
Finally we would like to thank to the anonymous referees, whose comments and suggestions have helped us to improve the
paper.

\pdfbookmark[1]{References}{ref}
\LastPageEnding

\end{document}